\documentclass[review]{elsarticle}
\usepackage{amssymb,amsmath,latexsym,amsfonts,bm,graphicx,mathrsfs,multirow,url,subfigure,algorithm,algpseudocode,booktabs,flafter}
\usepackage{color}
\usepackage{multicol}
\usepackage[center]{caption2}
\usepackage{mathdots,textcomp}
\usepackage{pifont,mathtools,epstopdf}
\usepackage{lineno,hyperref}
\modulolinenumbers[5]
\newtheorem{lemma}{Lemma}
\newtheorem{theorem}{Theorem}

\newtheorem{proposition}{Proposition}

\newtheorem{remark}{Remark}
\newtheorem{example}{Example}
\newenvironment{proof}{\noindent{\bf Proof.}}{\hfill\fbox{}\vspace*{1mm}}

\begin{document}

\begin{frontmatter}

\title{A fast two-level Strang splitting method for multi-dimensional spatial fractional Allen-Cahn equations with discrete maximum principle}

\author[Macau]{Yao-Yuan Cai}
\ead{yc07475@um.edu.mo}
\author[Foshan]{Zhi-Wei Fang}
\ead{fangzw@fosu.edu.cn}
\author[Chongqing]{Hao Chen}
\ead{hch@cqnu.edu.cn}
\author[Macau]{Hai-Wei Sun\corref{cor1}}
\cortext[cor1]{Corresponding author}
\ead{hsun@um.edu.mo}
\address[Macau]{Department of Mathematics, University of Macau, Macao, China.}
\address[Foshan]{School of Mathematics and Big Data, Foshan University, Foshan, 528000, China.}
\address[Chongqing]{College of Mathematics Science, Chongqing Normal University, Chongqing, 401331, China.}

\begin{abstract}
In this paper, we study the numerical solutions of the multi-dimensional spatial fractional Allen-Cahn equations. After semi-discretization for the spatial fractional Riesz derivative, a system of nonlinear ordinary differential equations with Toeplitz structure is obtained.  For the sake of reducing the computational complexity, a two-level Strang splitting method is proposed where the Toeplitz matrix in the system is split into the sum of a circulant matrix and a skew-circulant matrix. Therefore, the proposed method can be quickly implemented by the fast Fourier transform, substituting to calculate the expensive Toeplitz matrix exponential. Theoretically, the discrete maximum principle of our method is unconditionally preserved. Moreover, the analysis of error in the infinite norm with second-order accuracy is conducted in both time and space. Finally, numerical tests are given to corroborate our theoretical conclusions and the efficiency of the proposed method.
\end{abstract}


\end{frontmatter}

\section{Introduction}\label{sec1}
Over the last decades, because of the catholicity of anomalous diffusion phenomena in the real world and the nearly unavailable analytic solutions to most fractional diffusion equations, many investigators have developed efficient numerical methods to solve the fractional diffusion equations, especially considering the Riesz fractional derivative in space, which includes the discretization schemes \cite{Meerschaert-2004,ortigueira-2006,tian-2015,celik-2012,sousa-2015} and their fast computation \cite{huang-2020,huang-2022,huang-2022_2,huang-2021,lu-2021,pang-2016,lei-2010,pang-2012}.

In this paper, we focus on developing the numerical method for solving the following  high-dimensional spatial
fractional Allen-Cahn (SFAC) equations\cite{hou-2017}
\begin{align}\label{FAC}
	\left\{
	\begin{aligned}
		&u_t=\varepsilon^2\mathcal{L}_{x^d}^\alpha u+u-u^3, \qquad x \in \Omega , t \in (0,T],\\
		&u(x,0)=u^0(x), \qquad \qquad \;\;\; x\in \overline\Omega,\\
		&u\lvert_{\partial \Omega}=0, \qquad \qquad \qquad \qquad t \in (0,T],
	\end{aligned}\right.
\end{align}
where $\Omega=[a,b]^d$ is a  two-dimensional (2D) domain $(d=2)$ or three-dimensional (3D) domain $(d=3)$, $\varepsilon>0$ is an interfacial parameter, and $\mathcal{L}_{x^d}^\alpha$ denotes the $d$-dimensional Riesz fractional operator. For simplicity, we give the definition in one-dimensional case of order $\alpha_1\in(1,2)$ by
\begin{equation}\label{Riesz}
	\mathcal{L}_{x^1}^\alpha u
	=\mathcal{L}_{x^{(1)}}^{\alpha_1}u
	:=\frac{1}{-2\cos\frac{\alpha_1 \pi}{2}}
	\big({_a\mathcal{D}}^{\alpha_1}_{x^{(1)}}u
	+ {_{x^{(1)}}\mathcal{D}}_b^{\alpha_1}u
	\big),
\end{equation}
in which the left-side and right-side Riemann-Liouville fractional derivatives \cite{podlubny-1998} are defined as
\begin{align*}
	&{_a\mathcal{D}}_{x}^{\alpha}u
	=\frac{1}{\Gamma(2-\alpha)}\cdot\frac{\mathrm{d}^2}{\mathrm{d} {x}^2}\int_{a}^{x}\frac{u(\xi)}{({x}-\xi)^{\alpha-1}}\mathrm{d}\xi,\\
	&{_{x}\mathcal{D}}_b^{\alpha}u
	=\frac{1}{\Gamma(2-\alpha)}\cdot\frac{\mathrm{d}^2}{\mathrm{d} {x}^2}\int_{x}^{b}\frac{u(\xi)}{(\xi-{x})^{\alpha-1}}\mathrm{d}\xi.
\end{align*}
The 2D Riesz fractional derivative is defined as
$\mathcal{L}_{x^2}^\alpha u
=\mathcal{L}_{x^{(1)}}^{\alpha_1}u+\mathcal{L}_{x^{(2)}}^{\alpha_2}u$, and the 3D case is denoted by
$\mathcal{L}_{x^3}^\alpha u
=\mathcal{L}_{x^{(1)}}^{\alpha_1}u+\mathcal{L}_{x^{(2)}}^{\alpha_2}u+\mathcal{L}_{x^{(3)}}^{\alpha_3}u$,
where $ \alpha_1,\alpha_2,\alpha_3 \in (1,2)$ are the orders of the fractional derivative.

The Allen-Cahn equations have been widely used in many areas, such as two-phase incompressible fluids, mean curvature flows, dendritic growth dynamics, image inpainting, and segmentation\cite{lee_kim-2019,li-2012,du-2004,evans-1992,evans-1991,liu-2003}.
Meanwhile, the time-discrete formats of the SFAC equations have been investigated in recent years.
In 2017, Hou, Tang, and Yang \cite{hou-2017} presented a Crank-Nicolson scheme which is a second-order time discretization method and demonstrates the discrete maximum principle for the first time. However, this scheme is very computationally expensive in calculating nonlinear terms due to using the iteration algorithm, and it preserves the discrete maximum principle and energy decay conditionally.
In 2019, an operator splitting with alternating direction implicit (splitting-ADI) method\cite{he-2020} is proposed by He, Pan, and Hu. The splitting-ADI method is much faster than the previous scheme and extends to the fourth-order accuracy for time, but it fulfills the discrete maximum principle by restricting the time step.
Then, Du, Ju, Li, and Qiao provided a second-order exponential time differencing with the Runge-Kutta method (ETDRK2) \cite{du-2019}. The maximum principle is preserved by this method unconditionally. Nevertheless, this method is only suitable for periodic boundary conditions and the accuracy of the results is insufficient.
Recently, a dimensional splitting ETDRK2 method \cite{chen-2021-1} and a Strang splitting method\cite{chen-2021-2} were presented by Chen and Sun to solve the problem with Dirichlet boundary conditions efficiently. The discrete maximum principle has been proved theoretically. However, these methods take an amount of computing time as calculating the exponential of Toeplitz matrices by using the Gohberg-Semencul method.

In this paper, we concentrate on the fast computation of the discretization scheme of the fractional operator. A two-level Strang splitting method is proposed to avoid calculating the Toeplitz matrix exponential (TME), aiming at reducing the computational complexity and maintaining the discrete maximum principle. Noting the disadvantage in the Strang splitting approach\cite{chen-2021-2}, the linear term, where is a Toeplitz matrix from the space discretization, takes much CPU time to obtain the numerical solution. Significantly, any Toeplitz matrix can be expressed as a sum of a circulant matrix and a skew-circulant matrix \cite{ng-2003,chan-1993} which can be diagonalized efficiently. Thus, the Strang splitting method is utilized to the linear part once more, resulting in a two-level Strang splitting method. The benefit of the two-level Strang splitting algorithm is that it not only retains the second-order convergence in time but also greatly accelerates the calculation. The computational cost is significantly reduced to $\mathcal{O}(M\log M)$ by applying the fast Fourier transform (FFT) where $M$ denotes the spatial size in total. Besides, the discrete maximum principle, which states the absolute value of the entire solution is limited to 1 if the absolute value of initial and boundary data are bounded by 1, is still preserved unconditionally.

The remainder of the paper is laid out as follows. In Section 2, the fully discretized technique for approximating the SFAC problem will be shown, which includes a second-order finite difference in space and the two-level Strang splitting method in time. The proof of the discrete maximum principle of the proposed method is given in Section 3. The proposed method is studied under convergence analysis in time with norms of Banach space in Section 4. The numerical results are presented in Section 5 and some conclusions are given in Section 6.

\section{Fully discrete schemes of the SFAC equation}\label{FDM}
In the following, we introduce the numerical method to solve the multidimensional SFAC equations. The proposed method is based on a second-order finite difference approximation \cite{tian-2015} for space and a two-level Strang splitting method, which is second-order accurate in time. Here we consider the 2D case in detail, the 3D case is omitted for similarity.

\subsection{Second-order semi-discrete scheme}
We consider the finite difference approximation in space on uniform meshes. For any positive integers $m_\ell,~\ell=1,2,3$, the grid is divided as
\begin{equation*}
x^{(\ell)}_i=a+ih_\ell ,\quad h_\ell=\frac{b-a}{m_\ell}, \quad i=0,1,...,m_\ell.
\end{equation*}
In the following, we introduce the discretization in the view of one dimension, using $x$, $\alpha$, $m$, and $h$ substituting for $x^{(\ell)}$, $\alpha_\ell$, $m_\ell$, and $h_\ell$, respectively. According to the approximation in \cite{tian-2015}, the Riemann-Liouville fractional derivatives on the mesh are approximated as:
\begin{align*}
&_a\mathcal{D}_{x}^{\alpha}u(x,t)\big \lvert \lvert_{x=x_i}
=\frac{1}{h^{\alpha}}\sum_{k=0}^{i+1} \omega_{k}^{(\alpha)}u(x_{i-k+1},t)
 +\mathcal O(h^2),\\
&_x\mathcal{D}_b^{\alpha}u(x,t)\big\lvert_{x=x_i}
=\frac{1}{h^{\alpha}}\sum_{k=0}^{m-i+1} \omega_{k}^{(\alpha)}u(x_{i+k-1},t)
 +\mathcal O(h^2),
\end{align*}
where the coefficients $\omega_k^{(\alpha)}$ are composed of the alternating fractional binomial coefficient,
\begin{align}\label{omega_coe}
\left\{
\begin{aligned}
&\omega_0^{(\alpha)}=\frac{\alpha}{2}g_0^{(\alpha)}, \\
&\omega_{k+1}^{(\alpha)}
=\frac{\alpha}{2}g_{k+1}^{(\alpha)}
 +\frac{2-\alpha}{2}g_{k}^{(\alpha)}, \qquad k=0,1,2,...,
\end{aligned}
\right.
\end{align}
and
\begin{align}\label{g_coe}
\left\{
\begin{aligned}
&g_0^{(\alpha)}=1, \\
&g_{k+1}^{(\alpha)}
=\left(1-\frac{\alpha+1}{k+1}\right)g_{k}^{(\alpha)},
 \qquad k=0,1,2,....
\end{aligned}
\right.
\end{align}
The coefficients $\omega_{k}^{(\alpha)}$ have the following properties:
\begin{proposition}\label{omega_pro}
(see \cite{tian-2015}) Let $1<\alpha<2$ and $\omega_{k}^{(\alpha)}$ be defined in (\ref{omega_coe}). We have
\begin{align*}
\left\{
\begin{aligned}
&\omega_{0}^{(\alpha)}=\frac{\alpha}{2},\quad \omega_{1}^{(\alpha)}=\frac{2-\alpha-\alpha^2}{2} \leq 0,\quad \omega_{2}^{(\alpha)}=\frac{\alpha(\alpha^2+\alpha-4)}{4},\quad \omega_{0}^{(\alpha)}+\omega_{2}^{(\alpha)}>0,  \\
& 1 \geq \omega_{0}^{(\alpha)} \geq \omega_{3}^{(\alpha)} \geq \omega_{4}^{(\alpha)} \geq ... \geq 0,\\
&\sum_{k=0}^{\infty} \omega_{k}^{(\alpha)}= 0, \quad \sum_{k=0}^{m} \omega_{k}^{(\alpha)}< 0, \quad m \geq 2.
\end{aligned}
\right.
\end{align*}
\end{proposition}

Then, denoting $u_{i,j}(t)$ as an approximation of $u(x^{(1)}_i,x^{(2)}_j,t)$, we can give the semi-discrete second-order finite difference scheme of (\ref{FAC}) in 2D case:
\begin{align}\label{disFAC}
\begin{aligned}
\frac{\partial u_{i,j}(t)}{\partial t}=
&-\frac{\varepsilon^2}{2h_1^{\alpha_1}\cos\frac{\alpha_1\pi}{2}}
\bigg(\sum_{k_1=0}^{i+1}\omega_{k_1}^{(\alpha_1)}u_{i-k_1+1,j}(t)
+\sum_{k_1=0}^{m_1-i+1}\omega_{k_1}^{(\alpha_1)}u_{i+k_1-1,j}(t)\bigg)\\
&-\frac{\varepsilon^2}{2h_2^{\alpha_2}\cos\frac{\alpha_2\pi}{2}}
\bigg(\sum_{k_2=0}^{j+1}\omega_{k_2}^{(\alpha_2)}u_{i,j-k_2+1}(t)
+\sum_{k_2=0}^{m_2-j+1}\omega_{k_2}^{(\alpha_2)}u_{i,j+k_2-1}(t)\bigg)\\
&+u_{i,j}(t)-u_{i,j}(t)^3, \qquad 1\leq i\leq m_1-1, \quad 1\leq j\leq m_2-1,
\end{aligned}
\end{align}
with the initial boundary conditions:
\begin{equation}\label{condition}
	u_{0,j}(t)=u_{m_1,j}(t)=u_{i,0}(t)=u_{i,m_2}(t)=0, \qquad u_{i,j}(0)=u^0\big(x^{(1)}_i,x^{(2)}_j\big).
\end{equation}
To rewrite the numerical scheme in matrix form, we denote
$$
{\bf u}(t)
=[u_{1,1}(t),...,u_{m_1-1,1}(t),u_{1,2}(t),...,u_{m_1-1,2}(t),...,u_{m_1-1,m_2-1}(t)]^\intercal.
$$
Thus, we have
\begin{equation}\label{simpleFAC}
\frac{\partial {\bf u}(t)}{\partial t}=A{\bf u}(t)+f\big({\bf u}(t)\big),
\quad {\bf u}(0)={\bf u}^0,
\quad 0\leq t\leq T,
\end{equation}
with the nonlinear term $f\big({\bf u}(t)\big)={\bf u}(t)-{\bf u}(t)^{3}$ where ${\bf u}(t)^{3}$ means the $3$rd Hadamard power of ${\bf u}(t)$. The matrix A is presented as:
\begin{equation}\label{2DA}
A=I_{m_2} \otimes B_{\alpha_1}+ B_{\alpha_2} \otimes I_{m_1},
\end{equation}
where $I_{m_1}$ and $I_{m_2}$ are both identity matrices in each size and
\begin{equation}\label{B}
B_\alpha=-\frac{\varepsilon^2}{2h^{\alpha}\cos\frac{\alpha\pi}{2}}(D_{\alpha}+D_{\alpha}^\intercal),
\end{equation}
with
\begin{equation*}
D_{\alpha}=\left[
\begin{array}{cccccc}
\omega_1^{(\alpha)}&\omega_0^{(\alpha)}&0&\cdots&0&0\\
\omega_2^{(\alpha)}&\omega_1^{(\alpha)}&\omega_0^{(\alpha)}&0&\cdots&0\\
\vdots&\omega_2^{(\alpha)}&\omega_1^{(\alpha)}&\ddots&\ddots&\vdots\\
\vdots&\ddots&\ddots&\ddots&\ddots&0\\
\omega_{m-2}^{(\alpha)}&\ddots&\ddots&\ddots&\omega_1^{(\alpha)}& \omega_0^{(\alpha)}\\
\omega_{m-1}^{(\alpha)}&\omega_{m-2}^{(\alpha)}&\cdots&\cdots& \omega_2^{(\alpha)}&\omega_1^{(\alpha)} \\
\end{array}
\right].
\end{equation*}
Meanwhile, $B_\alpha$ is a symmetric Toeplitz matrix.

Similarly, in the 3D situation, the matrix A of the semi-discrete scheme can be exhibited as the following form:
\begin{equation}\label{3DA}
A=I_{m_3} \otimes I_{m_2} \otimes B_{\alpha_1}
 + I_{m_3} \otimes B_{\alpha_2} \otimes I_{m_1}
 + B_{\alpha_3} \otimes I_{m_2} \otimes I_{m_1}.
\end{equation}

\subsection{Two-level Strang splitting method}
To solve the semi-discrete system \eqref{simpleFAC}, we propose to deal with linear and nonlinear parts separately, which leads to the Strang splitting method \cite{strang-1968}. Precisely, the two subproblems are defined as:

Linear subproblem:
\begin{equation}\label{linear}
\frac{\partial {\bf v}_1(t)}{\partial t}=A{\bf v}_1(t), \qquad {\bf v}_1(0)={\bf v}_1^0.
\end{equation}

Nonlinear subproblem:
\begin{equation}\label{nonlinear}
\frac{\partial {\bf v}_2(t)}{\partial t}=f\big({\bf v}_2(t)\big)={\bf v}_2(t)-{\bf v}_2(t)^{3}, \qquad {\bf v}_2(0)={\bf v}_2^0.
\end{equation}
Denote that $\mathcal{J}_\tau^P$ is the linear solution operator and $\mathcal{J}_\tau^Q$ is the nonlinear solution operator. The solution of the linear subproblem \eqref{linear} is
\begin{equation}\label{linear_operator}
{\bf v}_1(\tau)=e^{\tau A}{\bf v}_1^0:=\mathcal{J}_\tau ^P{\bf v}_1^0.
\end{equation}
As for the second part (\ref{nonlinear}) which is the classic Bernoulli differential equation \cite{hairer-1993}, the solution is given below:
\begin{equation}\label{nonlinear_operator}
{\bf v}_2(\tau)=\frac{{\bf v}_2^0}{\sqrt{({\bf v}_2^0)^{2}
		        +\big[\mathbf{1}_{m-1}-({\bf v}_2^0)^{2}\big]e^{-2\tau}}}
	           :=\mathcal{J}_\tau^Q{\bf v}_2^0.
\end{equation}
where  ${\bf1}_{m-1}=[1,1,...,1]^\intercal\in\mathbb{R}^{m-1}$.

Then, let $\tau$ be the size of time step and define $t_n=\tau n \leq T$ as the temporal partition for $n=0,1,...$.
Denote that ${\bf u}^{n}\approx {\bf u}(t_n)$ as the numerical approximation at a certain point $t_n$ by using the Strang splitting method. Note that the initial boundary conditions are remained, which means that ${\bf u}^{0}={\bf u}(t_0)$. Hence we have the following temporal second-order scheme and its recursive form:
\begin{equation}\label{strang}
{\bf u}^{n+1}=\mathcal{M}_\tau {\bf u}^n
=\mathcal{J}_\frac{\tau}{2}^Q
 \mathcal{J}_\tau^P\mathcal{J}_\frac{\tau}{2}^Q{\bf u}^n
=(\mathcal{J}_\frac{\tau}{2}^Q
 \mathcal{J}_\tau^P\mathcal{J}_\frac{\tau}{2}^Q)^{n+1}{\bf u}^0.
\end{equation}

Aiming at proposing a fast algorithm, it is necessary to give an efficient technique approximating the linear subproblem. In the literature \cite{ng-2003,chan-1993}, we could exploit a split of the matrix $B_\alpha=C_\alpha+S_\alpha$, where
\begin{equation}\label{B}
B_\alpha
=\left[
\begin{array}{ccccc}
b_0&b_{-1}&b_{-2}&\cdots&b_{2-m}\\
b_1&b_0&b_{-1}&\cdots&b_{3-m}\\
b_2&b_1&b_0&\cdots&b_{4-m}\\
\vdots&\ddots&\ddots&\ddots&\vdots\\
b_{m-2}&b_{m-3}&\cdots&\cdots& b_{0} \\
\end{array}
\right],
\end{equation}
i.e., $b_{ij}=b_{i-j}$ and $B_\alpha$ is constant along its diagonals, and the $k$th diagonal entries of $C_\alpha=[c_{i-j}]$ and $S_\alpha=[s_{i-j}]$ is defined by
\begin{equation}\label{C}
c_k=
\frac{1}{2}
\begin{cases}
b_k+b_{k-m+1},&k=1,2\ldots,m-2,\\
b_0,&k=0,\\
b_k+b_{k+m-1},&k=-1,-2\ldots,2-m,
\end{cases}
\end{equation}
and
\begin{equation}\label{S}
s_k=
\frac{1}{2}\begin{cases}
b_k-b_{k-m+1},&k=1,2\ldots,m-2,\\
b_0,&k=0,\\
b_k-b_{k+m-1},&k=-1,-2\ldots,2-m.
\end{cases}
\end{equation}

We rewrite the linear subproblem (\ref{linear}) to the following form:
\begin{equation}\label{newlinear}
\frac{\partial {\bf v}_1(t)}{\partial t}
=(S_\alpha+C_\alpha){\bf v}_1(t), \qquad {\bf v}_1(0)={\bf v}_1^0.
\end{equation}
where $C_\alpha=I_{m_2} \otimes C_{\alpha_1}+ C_{\alpha_2} \otimes I_{m_1}$, $S_\alpha=I_{m_2} \otimes S_{\alpha_1}+ S_{\alpha_2} \otimes I_{m_1}$ in the 2D case
and $ C_\alpha=I_{m_3} \otimes I_{m_2} \otimes C_{\alpha_1}
+ I_{m_3} \otimes C_{\alpha_2} \otimes I_{m_1}
+ C_{\alpha_3} \otimes I_{m_2} \otimes I_{m_1}
$,
$ S_\alpha=I_{m_3} \otimes I_{m_2} \otimes S_{\alpha_1}
+ I_{m_3} \otimes S_{\alpha_2} \otimes I_{m_1}
+ S_{\alpha_3} \otimes I_{m_2} \otimes I_{m_1}$ in the 3D case. Denoting that ${\bf \hat{v}}_1^{n+1} \approx {\bf v}_1^{n+1}$ is the numerical approximation to the solution of the linear subproblem (\ref{newlinear}), a second-level Strang splitting is employed for speeding up the numerical solution by the recursion. Here is the expression of the numerical approximation form:
\begin{equation}\label{hatv}
{\bf \hat{v}}_1^{n+1}
=\widehat{\mathcal{J}}_\tau^P{\bf \hat{v}}_1^n
=\left(e^{\frac{\tau}{2}C_\alpha}e^{\tau S_\alpha}
 e^{\frac{\tau }{2}C_\alpha}\right){\bf \hat{v}}_1^n.
\end{equation}

Consequently, an improved algorithm for calculating the equation (\ref{strang}),  which is called a two-level Strang splitting method, is presented below:
\begin{equation}\label{finalstrang}
{\bf \hat{u}}^{n+1}=\widehat{\mathcal{M}}_\tau {\bf \hat{u}}^n
=\mathcal{J}_\frac{\tau}{2}^Q
 \widehat{\mathcal{J}}_\tau^P\mathcal{J}_\frac{\tau}{2}^Q{\bf \hat{u}}^n
=(\mathcal{J}_\frac{\tau}{2}^Q
 \widehat{\mathcal{J}}_\tau^P\mathcal{J}_\frac{\tau}{2}^Q)^{n+1}{\bf u}^0.
\end{equation}
\begin{remark}
We remark that the circulant matrix and skew-circulant matrix could be diagonalized by use of the Fourier matrix. It is well-known \cite{chan-1996} that any circulant matrix can be diagonalized by the Fourier matrix $F$, i.e.,
$$
C_\alpha=F^*\Lambda_{C_\alpha}F,
$$
and a skew-circulant matrix has the spectral decomposition:
$$
S_\alpha=\Psi^*F^*\Lambda_{S_\alpha}F\Psi,
$$
where $F \in \mathbb{R}^{m-1 \times m-1}$ is a Fourier matrix and the entries of $F$ are given by $(F)_{jk}=\frac{1}{\sqrt{m-1}}e^{\frac{2\pi {\bf{i}}jk}{m-1}}$ for $0\leq j,k \leq m-2$ and $\Psi=\text{diag}\{1,e^{-\frac{{\bf{i}}\pi}{m-1}},\cdots,e^{-\frac{{\bf{i}}\pi(m-2)}{m-1}}\}$ is a diagonal matrix. Note that $F$ is unitary, and its first column is $\frac{{\bf1}_{m-1}}{\sqrt {m-1}}$. Then:
$$
FC_\alpha {\bf e}_1=\frac{1}{\sqrt{m-1}}\Lambda_{C_\alpha}{\bf1}_{m-1}, \qquad
F\Psi S_\alpha {\bf e}_1=\frac{1}{\sqrt{m-1}}\Lambda_{S_\alpha}{\bf1}_{m-1},
$$
where ${\bf e}_1$ is the first unit vector. Therefore, the proposed linear subproblem \eqref{newlinear} can be carried out efficiently by the FFTs.
\end{remark}

\section{Discrete maximum principle}\label{Discrete maximum principle}
In this section, we prove one of the particularly important physical properties, discrete maximum principle, which means the solution bounded by the initial and boundary data. Before proving, a definition will be presented below for the following lemmas and proofs.
\begin{lemma}
(see\cite{powers-1975}) Without loss of generality, if a matrix $L$ is a diagonally dominant matrix with the main nonpositive diagonal entries, then
$$
\|e^{\tau L}\|_{\infty} \leq 1, \qquad \tau \geq 0.
$$
\end{lemma}

\begin{lemma}
Let the circulant matrix $C_\alpha$ and the skew circulant matrix $S_\alpha$ be defined by \eqref{C} and \eqref{S}, respectively. Then,  $C_\alpha$ and $S_\alpha$ are diagonally dominant matrices with negative diagonal entries.
\end{lemma}
\begin{proof}
First, we consider the circulant matrix $C_\alpha$. By the definition \eqref{C}, and the properties of $\omega_k^{(\alpha)}$ in Proposition \ref{omega_pro}. We know that the entries of $C_\alpha$ are positive except the negative main diagonals. Meanwhile, any circulant matrix has the same row sums. The row sums of $C_\alpha$ is derived as
\begin{equation*}
\sum_{k=0}^{m-2}c_k
=\frac{1}{2}\sum_{k=2-m}^{m-2} b_k
=-\frac{\varepsilon^2}{2h^{\alpha}\cos\frac{\alpha\pi}{2}}\sum_{k=0}^{m-1}\omega_{k}^{(\alpha)}\leq 0.
\end{equation*}
We deduce that the matrix $C_\alpha$ meets the lemma.
	
Afterward, the skew-circulant matrix $S_\alpha$ has negative main diagonals. Notice that the sums of the  absolute value of each row are the same. Thus, we have
\begin{equation*}
\sum_{k=1}^{m-2}\big\lvert s_k\big \lvert
\leq\frac{1}{2}\sum_{k=2-m,k\neq0}^{m-2}\big\lvert b_k\big \lvert
=-\frac{\varepsilon^2}{2h^{\alpha}\cos\frac{\alpha\pi}{2}}
\sum_{k=0,k\neq1}^{m-1}\omega_{k}^{(\alpha)}\leq\big\lvert s_0\big\lvert,
\end{equation*}
which finishes the proof.
\end{proof}

Accordingly, we conclude that $\|e^{\tau C_\alpha }\|_{\infty} \leq 1$ and $\|e^{\tau S_\alpha}\|_{\infty} \leq 1$ for any $\tau \geq 0$ in (\ref{newlinear}). From the above lemmas, we have the following lemma.
\begin{lemma}\label{lemma:linear}
Define $\widehat{\mathcal{J}}_\tau^P$ as the solution operator of the linear subproblem in (\ref{hatv}), then for any time step $\tau\geq0$ and ${\bf v}\in \mathbb{R}^{m}$ with $\|{\bf v}\|_{\infty}\leq 1$, we have the following inequality:
\begin{equation}
\|\widehat{\mathcal{J}}_\tau^P{\bf v}\|_{\infty} \leq 1.
\end{equation}
\end{lemma}

Next, we consider the nonlinear subproblem for investigating the discrete maximum principle.
\begin{lemma}\label{lemma:nonlinear}
(See \cite{chen-2021-1}) Define $\mathcal{J}_\tau^Q$ as the solution operator of the nonlinear subproblem in (\ref{nonlinear_operator}). It holds that:
$$
\|\mathcal{J}_\tau^Q{\bf v}\|_{\infty} \leq 1,
$$
for any time step $\tau\geq 0$ and ${\bf v}\in \mathbb{R}^{m}$ with $\|{\bf v}\|_{\infty}\leq 1$.
\end{lemma}

Now we start to prove the discrete maximum principle.
\begin{theorem} 
Suppose that ${\bf u}^0$ is an initial input data that satisfies $\|{\bf u}^0\|_{\infty} \leq 1$. The proposed method \eqref{finalstrang} preserves the discrete maximum principle in any arbitrary time step size $\tau>0$, i.e.,
$$
\|{\bf u}^n\|_{\infty}\leq 1, \quad \forall n \geq 0.
$$
\end{theorem}
\begin{proof}
By mathematical induction, we assume $\|{\bf u}^k\|_{\infty} \leq 1$. Then, it is to show the result correct when $n=k+1$. According to the two-level Strang splitting scheme, we obtain
\begin{equation}
\|{\bf u}^{k+1} \|_{\infty}
=\big\|\mathcal{J}_\frac{\tau}{2}^Q\widehat{\mathcal{J}}_\tau^P
\mathcal{J}_\frac{\tau}{2}^Q{\bf u}^k \big\|_{\infty}
\leq 1.
\end{equation}
The theorem is concluded by use of Lemma \ref{lemma:linear} and Lemma \ref{lemma:nonlinear}.
\end{proof}

\section{Convergence analysis}
In this section, we consider the convergence of the proposed method from the SFAC equations. The second-order convergence in the spatial discretization is verified in articles \cite{chen-2021-1,chen-2021-2}. We concentrate on the temporal convergence order. Some lemmas will be presented for the convergence of the two-level Strang splitting technique in the following proof.
\begin{lemma}\label{lemma:first_level}
(See \cite{hansen-2012}) For the semilinear equation \eqref{simpleFAC}, assume that:
	
(\romannumeral1) $\|e^{\tau A}\|_\infty\leq 1$ for any $\tau \geq 0$.
	
(\romannumeral2) The nonlinear function f is two times continuously Fr\'{e}chet differentiable for all ${\bf u}$.
	
(\romannumeral3) The exact solution ${\bf u}$ of the semilinear equation \eqref{simpleFAC} belongs to $C^2([0,T])$.\\
Then the Strang splitting method \eqref{strang} satisfies the second-order consistency, i.e.,
\begin{equation}\label{consist}
\|{\bf u}(t+\tau)-\mathcal{M}_\tau{\bf u}(t)\|_\infty \leq g_1\tau^3,
\end{equation}
for a positive constant $g_1$.
\end{lemma}

\begin{lemma}\label{lemma:linear_v}
For the linear equation
\begin{equation}\label{test}
\frac{\partial {\bf w}_1(\tau)}{\partial \tau}
=(S+C){\bf w}_1(\tau), \qquad {\bf w}_1(0)={\bf w}_1^0.
\end{equation}
The third-order expansion of the exact solution is then given by
{\footnotesize \begin{align}\label{eq:linear_ex}
\begin{aligned}		
{\bf w}_1(\tau)
=&e^{S\tau}{\bf w}_1^0
 +\int_{0}^{\tau}e^{S(\tau-\delta_1)}Ce^{S\delta_1}{\bf w}_1^0d\delta_1
 +\int_{0}^{\tau}\int_{0}^{\delta_1}e^{S(\tau-\delta_1)}Ce^{S(\delta_1-\delta_2)}Ce^{S\delta_2}{\bf w}_1^0d\delta_2d\delta_1\\
&+\int_{0}^{\tau}\int_{0}^{\delta_1}\int_{0}^{\delta_2}e^{S(\tau-\delta_1)}Ce^{S(\delta_1-\delta_2)}Ce^{S(\delta_2-\delta_3)}Ce^{S\delta_3}{\bf w}_1(\delta_3)
d\delta_3d\delta_2d\delta_1,
\end{aligned}
\end{align}}
where $\tau \geq \delta_1 \geq \delta_2 \geq \delta_3 \geq 0$.
	
\begin{proof}
The lemma is proved by applying the variation-of-constant formula three times.
\end{proof}
\end{lemma}

\begin{lemma}\label{lemma:CSC_v}
(See \cite{einkemmer-2014}) The splitting operator (\ref{hatv}) has the expansion as below:
\begin{equation}\label{eq:linear_aprox}
\widehat{\mathcal{J}}_\tau^P{\bf w}_1^0
=e^{S\tau}{\bf w}_1^0
+\frac{\tau}{2}\big\{C,e^{S\tau}\big\}{\bf w}_1^0
+\frac{\tau^2}{8}\Big\{C,\big\{C,e^{S\tau}\big\}\Big\}{\bf w}_1^0
+R_3{\bf w}_1^0
\end{equation}
where
$$
R_3=\frac{\tau^3}{16}\int_{0}^{1}(1-\theta)^2
\bigg\{C,\Big\{C,\big\{C,e^{\frac{\tau\theta}{2}C}e^{\tau S}e^{\frac{\tau\theta}{2}C}\big\}\Big\}\bigg\}d\theta,
$$
$\tau\geq\theta\geq 0$, and $\{C,S\}=CS+SC$.
\end{lemma}

\begin{lemma}\label{the:second_level}
Considering the linear subproblem \eqref{test}, suppose that
\begin{align}
\label{eq1}
\Big\|e^{S(\tau-\delta_1)}\big[[C,S],S\big] e^{S\delta_1}
{\bf w}_1^0\Big\|_\infty &\leq g;\\
\label{eq2} 
\Big\|e^{S(\tau-\delta_1)}\big[[C,S],S\big]
e^{S(\delta_1-\delta_2)}Ce^{S\delta_2}{\bf w}_1^0\Big\|_\infty &\leq g;\\
\label{eq3}
\Big\|e^{S(\tau-\delta_1)}Ce^{S(\delta_1-\delta_2)}
\big[[C,S],S\big]e^{S\delta_2}{\bf w}_1^0\Big\|_\infty &\leq g;\\
\label{eq4}
\bigg\|\Big\{C,\big\{C,\{C,e^{\frac{\phi}{2}C}e^{\tau S}
e^{\frac{\phi}{2}C}\}\big\}\Big\}{\bf w}_1^0\bigg\|_\infty &\leq g;\\
\label{eq5}
\Big\|e^{S(\tau-\delta_1)}Ce^{S(\delta_1-\delta_2)}C
e^{S(\delta_2-\delta_3)}Ce^{S\delta_3}{\bf w}_1(\delta_3)\Big\|_\infty &\leq g,
\end{align}
where $[C,S]=CS-SC$ and a positive constant $g$. For sufficiently small time step $\tau\geq\delta_1\geq\delta_2\geq\delta_3\geq0$, $\phi \in[0,\tau]$, and the input ${\bf w}_1^0$, the splitting operator \eqref{eq:linear_aprox} is consistent with order 2, i.e.,
\begin{equation*}
\|{\bf w}_1(\tau)-\widehat{\mathcal{J}}_\tau^P{\bf w}_1^0\|_\infty
\leq g_2\tau^3.
\end{equation*}
\begin{proof}
Similarly to the reference \cite{einkemmer-2014}, we show the difference between ${\bf w}_1(\tau)$ in Lemma \ref{lemma:linear_v} and $\widehat{\mathcal{J}}_\tau^P{\bf w}_1^0$ in Lemma \ref{lemma:CSC_v}. It is to compare the same order terms on the right sides of \eqref{eq:linear_ex} and \eqref{eq:linear_aprox}. Here we start from second-order items as first-order items are the same.
		
\textbf{Second-order items.} Setting $p_1(\delta_1)=e^{S(\tau-\delta_1)}Ce^{S\delta_1}{\bf w}_1^0$, we get by using the trapezoidal rule \cite{jahnke-2000}
$$
\Big\|\frac{\tau}{2}\big[p_1(0)+p_1(\tau)\big]-\int_{0}^{\tau}p_1(\theta)d\theta\Big\|_\infty
=\frac{\tau^3}{12}\left\|p_1''(\tilde{\delta})\right\|_\infty,
$$
where $\tilde{\delta} \in [0,\tau]$. The residual term is bounded by the assumption (\ref{eq1}) for any $\delta_1$.
		
\textbf{Third-order items.} For
$p_2(\delta_1,\delta_2)
=e^{S(\tau-\delta_1)}Ce^{S(\delta_1-\delta_2)}Ce^{S\delta_2}{\bf w}_1^0$, we adopt a trapezoidal rule and mean value theorem
\begin{align*}
&\Big\|\frac{\tau^2}{8}\big\{C,\{C,e^{S\tau}\}\big\}{\bf w}_1^0
 -\int_{0}^{\tau}\int_{0}^{\delta_1}e^{S(\tau-\delta_1)}C
e^{S(\delta_1-\delta_2)}Ce^{S\delta_2}
{\bf w}_1^0d\delta_2d\delta_1\Big\|_\infty\\
=&\Big\|\frac{\tau^2}{8} \big[p_2(\tau,\tau)+2p_2(\tau,0)+p_2(0,0)\big]
-\int_{0}^{\tau}\int_{0}^{\delta_1}p_2(\delta_1,\delta_2)d\delta_2d\delta_1
\Big\|_\infty\\
=&\frac{\tau^4}{24}\Big\|
\frac{\partial^2 p_2}{\partial \delta_2^2}(\theta,\zeta_2)
+\frac{1}{2}\Big(\frac{\partial^2 p_2}{\partial \delta_1^2}(\tau,\zeta_1)
+\frac{\partial^2 p_2}{\partial \delta_1^2}(\zeta_1,0)\Big)
\Big\|_\infty
\end{align*}
where $0\leq\zeta_2\leq\theta\leq\tau$ and $0\leq \zeta_1\leq \tau$, with
{\footnotesize
\begin{align*}
\frac{\partial^2 p_2}{\partial \delta_1^2}
=e^{S(\tau-\delta_1)}\big[[C,S],S\big]
e^{S(\delta_1-\delta_2)}Ce^{S\delta_2}{\bf w}_1^0,~
\frac{\partial^2 p_2}{\partial \delta_2^2}
=e^{S(\tau-\delta_1)}Ce^{S(\delta_1-\delta_2)}\big[[C,S],S\big]e^{S\delta_2}{\bf w}_1^0.
\end{align*}}
This remainder terms are bounded if the suppositions \eqref{eq2} and \eqref{eq3} holds.
		
\textbf{Fourth-order items.} For bounding the residual terms in the expansion of the exact solution and the approximate solution, we need the presumptions \eqref{eq4} and \eqref{eq5}.
\end{proof}
\end{lemma}
From the Section \ref{Discrete maximum principle}, $\|e^{ \tau C_{\alpha}}\|_\infty \leq 1$ and $\|e^{\tau S_{\alpha}}\|_\infty \leq 1$. Furthermore, in order to maintain the stability, the ratio of space step size and time step size should be controlled, for example, $\frac{\tau}{h^2}\leq g$, where $g$ is an positive number. Therefore, there exist an positive number $g_3$ and we derive the following inequality in each $d-$dimensional cases, i.e.,
\begin{equation}
\|\mathcal{J}_\tau^P{\bf w}_1^0
-\widehat{\mathcal{J}}_\tau^P{\bf w}_1^0\|_\infty
\leq g_3\tau^3.
\end{equation}

\begin{lemma}
Let ${\bf v},{\bf w}\in\mathbb{R}^m$, then we have
\begin{equation*}
\|\mathcal{M}_\tau{\bf v}-\mathcal{M}_\tau{\bf w}\|_\infty
\leq e^\tau\|{\bf v}-{\bf w}\|_\infty,
\qquad
\|\widehat{\mathcal{M}}_\tau{\bf v}-\widehat{\mathcal{M}}_\tau{\bf w}\|_\infty
\leq e^\tau\|{\bf v}-{\bf w}\|_\infty.
\end{equation*}
\begin{proof}
The left inequality is proved by \cite{chen-2021-2}. Similarly, the right inequality is proved as follows:
\begin{align*}
\|\widehat{\mathcal{M}}_\tau{\bf v}-\widehat{\mathcal{M}}_\tau{\bf w}\|_\infty
&= \big\|\mathcal{J}_\frac{\tau}{2}^Q\widehat{\mathcal{J}}_\tau^P
\mathcal{J}_\frac{\tau}{2}^Q{\bf v}
-\mathcal{J}_\frac{\tau}{2}^Q\widehat{\mathcal{J}}_\tau^P
\mathcal{J}_\frac{\tau}{2}^Q{\bf w}\big\|_{\infty}\\
&\leq e^{\frac{\tau}{2}}
\big\|\widehat{\mathcal{J}}_\tau^P \mathcal{J}_\frac{\tau}{2}^Q{\bf v}
-\widehat{\mathcal{J}}_\tau^P\mathcal{J}_\frac{\tau}{2}^Q
{\bf w}\big\|_{\infty}\\
&\leq e^{\frac{\tau}{2}}
\big\| \mathcal{J}_\frac{\tau}{2}^Q{\bf v}
-\mathcal{J}_\frac{\tau}{2}^Q{\bf w}\big\|_{\infty}\\
&\leq e^{\tau} \big\|{\bf v}-{\bf w}\big\|_{\infty}.
\end{align*}
Then we complete the proof.
\end{proof}
\end{lemma}
Therefore we derive the following lemma:
\begin{lemma}\label{newconsistency}
If the two-level Strang splitting scheme satisfies these assumptions (\romannumeral1), (\romannumeral2), and (\romannumeral3) in Lemma \ref{lemma:first_level}, and obeys \eqref{eq1}, \eqref{eq2}, \eqref{eq3}, \eqref{eq4}, and \eqref{eq5} in Lemma \ref{the:second_level} when $\tau \geq 0$. Moreover, for an initial value ${\bf u}^0\in \overline\Omega$, we assume
$\|{\bf u}^0\|_\infty\leq 1$, then we have the following inequality:
\begin{equation*}
\|\mathcal{M}_\tau {\bf u}^0
-\widehat{\mathcal{M}}_\tau {\bf u}^0\|_\infty
\leq g_3e^{\frac{\tau}{2}}\tau^3,
\end{equation*}
where $g_3$ are positive numbers.
	
\begin{proof}
From Lemma \ref{the:second_level} and Lemma \ref{newconsistency}, we derive that
\begin{align*}
\|\mathcal{M}_\tau{\bf u}^0-\widehat{\mathcal{M}}_\tau{\bf u}^0\|_\infty
&= \big\|\mathcal{J}_\frac{\tau}{2}^Q\mathcal{J}_\tau^P
\mathcal{J}_\frac{\tau}{2}^Q{\bf u}^0
-\mathcal{J}_\frac{\tau}{2}^Q\widehat{\mathcal{J}}_\tau^P
\mathcal{J}_\frac{\tau}{2}^Q{\bf u}^0\big\|_{\infty}\\
&\leq e^{\frac{\tau}{2}}
\big\|\mathcal{J}_\tau^P \mathcal{J}_\frac{\tau}{2}^Q{\bf u}^0
-\widehat{\mathcal{J}}_\tau^P\mathcal{J}_\frac{\tau}{2}^Q
{\bf u}^0\big\|_{\infty}\\
&\leq e^{\frac{\tau}{2}}
\big\| \mathcal{J}_\tau^P{\bf v}^0
-\widehat{\mathcal{J}}_\tau^P{\bf v}^0\big\|_{\infty}\\
&\leq g_3e^{\frac{\tau}{2}}\tau^3.
\end{align*}
where ${\bf v}^0=\mathcal{J}_\frac{\tau}{2}^Q {\bf u}^0$. Then we complete the proof.
\end{proof}
\end{lemma}

\begin{theorem}
Let $\hat{\bf u}^n$ and $u_{\bf x}(t_n)$ be the answer of the fully discretized scheme \eqref{finalstrang} and the exact solution of the SFAC equation on the uniform mesh at the moment $t_n=n\tau\leq T$. Assume that $u(t_n)\in C^{5,2}(\overline{\Omega};[0,T])$ and the initial value $\|{\bf u}^0\|_\infty \leq 1$, then the fully discretization scheme \eqref{finalstrang} takes the second-order convergence in both time and space; i.e.,
\begin{equation*}
\|u_{\bf x}(t_n)-\hat{\bf u}^n\|_\infty\leq g_5(h^2+\tau^2).
\end{equation*}
\begin{proof}
Suppose ${\bf u}^n$ is the Strang splitting scheme \eqref{strang} and ${\bf u}(t_n)$ is the exact solution of the semilinear equation \eqref{simpleFAC}. We derive by utilizing the triangle inequality:
\begin{equation}
\|u_{\bf x}(t_n)-\hat{\bf u}^n\|_\infty
\leq\|\hat{\bf u}^n-{\bf u}^n\|_\infty
+\|{\bf u}^n-{\bf u}(t_n)\|_\infty
+\|{\bf u}(t_n)-u_{\bf x}(t_n)\|_\infty.
\end{equation}
First we consider the temporal convergence. Denote $\mathcal E_\tau$ as the analytic solution operator of the semilinear equation \eqref{simpleFAC} in each time step, so we derive the global error by employing the Lady Windermere's fan \cite{hairer-2006}:
\begin{align}
\label{1}
{\bf u}^n-{\bf u}(t_n)
&={\mathcal{M}}_\tau^n{\bf u}^0-\mathcal E_\tau^n{\bf u}^0
=\sum_{l=1}^{n}\big({\mathcal{M}}_\tau^{n-l}\big)
({\mathcal{M}}_\tau\mathcal E_\tau^{l-1}-\mathcal E_\tau^{l})({\bf u}^0),\\
\label{2}
\hat{\bf u}^n-{\bf u}^n
&=\widehat{\mathcal{M}}_\tau^n {\bf u}^0
 -{\mathcal{M}}_\tau^n{\bf u}^0
=\sum_{l=1}^{n}\big(\widehat{\mathcal{M}}_\tau^{n-l}\big)
 (\widehat{\mathcal{M}}_\tau{\mathcal{M}}_\tau^{l-1}
 -{\mathcal{M}}_\tau^l)({\bf u}^0).
\end{align}
Taking the infinite-norm between two-sides of equations (\ref{1},\ref{2}), we have
\begin{equation*}
\begin{aligned}
\|{\bf u}^n-{\bf u}(t_n)\|_\infty
&\leq\sum_{l=1}^{n}
\|\big({\mathcal{M}}_\tau^{n-l}\big)({\mathcal{M}}_\tau\mathcal E_\tau^{l-1}
-\mathcal E_\tau^{l})({\bf u}^0)\|_\infty\\
&\leq\sum_{l=1}^{n}e^{(n-l)\tau}
\|({\mathcal{M}}_\tau\mathcal E_\tau^{l-1}
-\mathcal E_\tau^{l}){\bf u}^0\|_\infty,
\end{aligned}
\end{equation*}
and
\begin{equation*}
\begin{aligned}
\|\hat{\bf u}^n-{\bf u}^n\|_\infty
&\leq\sum_{l=1}^{n}
\|\big(\widehat{\mathcal{M}}_\tau^{n-l}\big)
(\widehat{\mathcal{M}}_\tau{\mathcal{M}}_\tau^{l-1}
-{\mathcal{M}}_\tau^l)({\bf u}^0)\|_\infty\\
&\leq\sum_{l=1}^{n}e^{(n-l)\tau}
\|(\widehat{\mathcal{M}}_\tau{\mathcal{M}}_\tau^{l-1}
-{\mathcal{M}}_\tau^l){\bf u}^0\|_\infty.
\end{aligned}
\end{equation*}
According Lemma \ref{lemma:first_level} and Lemma \ref{newconsistency}, we have:
\begin{align*}
\Big\|\big({\mathcal{M}}_\tau\mathcal E_\tau^{l-1}-\mathcal E_\tau^{l}\big){\bf u}^0\Big\|_\infty
&= \Big\|{\mathcal{M}}_\tau {\bf u}(t_{l-1})-\mathcal E_\tau {\bf u}(t_{l-1})\Big\|_\infty
 \leq  g_1\tau^3,\\
\Big\|(\widehat{\mathcal{M}}_\tau{\mathcal{M}}_\tau^{l-1}
-{\mathcal{M}}_\tau^l){\bf u}^0\Big\|_\infty
&= \Big\|\widehat{\mathcal{M}}_\tau {\bf u}^{l-1}
-{\mathcal{M}}_\tau {\bf u}^{l-1}\Big\|_\infty
\leq  g_3e^{\frac{\tau}{2}}\tau^3.
\end{align*}
Therefore, we get the following form:
\begin{equation*}
\|{\bf u}^n-{\bf u}(t_n)\|_\infty
\leq \sum_{l=1}^{n}e^{(n-l)\tau}g_1\tau^3
\leq e^{t_n}g_1\tau^2,
\end{equation*}
\begin{equation*}
\|\hat{\bf u}^n-{\bf u}^n\|_\infty
\leq \sum_{l=1}^{n}e^{(n-l)\tau}g_3e^{\frac{\tau}{2}}\tau^3
\leq e^{\frac{3t_n}{2}}g_3\tau^2
\end{equation*}
as $e^{\tau}-1 \geq \tau$ for $\tau \geq 0$ in the last step.
		
Then we consider the spatial convergence, which is accord with the third part of the right-hand side. In the light of \cite{chen-2021-1,chen-2021-2}, if the exact solution $u(t_n)\in C^5(\overline{\Omega})$, it follows that
\begin{equation}
\frac{\partial u_{\bf x}(t)}{\partial t}
=Au_{\bf x}(t_n)+f\big(u_{\bf x}(t_n)\big)+{\bf R}_h,
\end{equation}
where ${\bf R}_h$ is the local truncation error and the inequality of its infinite norm is
$\|{\bf R}_h\|_\infty \leq g_4\varepsilon^2 h^2$, in which $g_4$ is a positive constant. Then we have
		
\begin{equation*}
\|{\bf u}(t_n)-u_{\bf x}(t_n)\|_\infty
\leq g_4\varepsilon^2 h^2e^{2t_n}.
\end{equation*}
The convergence order is proved with
\begin{equation*}
g_5=\max\{e^{t_n}g_1+e^{\frac{3t_n}{2}}g_3, g_4\varepsilon^2 e^{2t_n}\}.
\end{equation*}
\end{proof}
\end{theorem}

\section{Numerical results} \label{numerical results}
In this section, we give some examples to demonstrate the effectiveness of the proposed method given in (\ref{finalstrang}), and use some graphs to show that the resulting numerical solution satisfies some physical properties, such as discrete maximum principle and energy dissipation. All numerical experiments are run in MATLAB R2019a on a laptop with the configuration: Intel(R) Core(TM) i7-10750H CPU @ 2.60GHz and 32 GB RAM.

For highlighting the advantage of our proposed scheme, two other methods are performed below in 2D form and compared. 
The first method was developed in \cite{chen-2021-2}, which is the Strang splitting scheme \eqref{strang}. We focus on the fast computation of TME, where it exists in the linear solution operator. Recently, some investigators applied the Krylov subspace methods to the matrix exponential \cite{lee-2010,lee-2012,moret-2004,zhang-2015,zhang-2020,zhang-2020_2}, here we use the shift-invert Lanczos method\cite{pang-2011,pang-2014} to calculate the TME.

The second method was proposed by \cite{he-2020}, the splitting-ADI scheme:
\begin{small}
\begin{equation*}
{\bf u}^{n+1}=\mathcal{J}_\frac{\tau}{2}^Q
\bigg\{\Big[(I_{m_2}-\frac{\tau}{2}B_{\alpha_2})^{-1}(I_{m_2}+\frac{\tau}{2}B_{\alpha_2})\Big]
\otimes
\Big[(I_{m_1}-\frac{\tau}{2}B_{\alpha_1})^{-1}(I_{m_1}+\frac{\tau}{2}B_{\alpha_1})\Big]\bigg\}
\mathcal{J}_\frac{\tau}{2}^Q{\bf u}^n.
\end{equation*}
\end{small}
These two methods mentioned above all encounter the problem of how to calculate the matrix-vector multiplication of the inverse of the Toeplitz matrix. Therefore, the Gohberg-Semencul-type formula \cite{gohberg-1992} to compute the matrix-vector production of the inverse of Toeplitz matrix.
In addition, the shift-invert Lanczos algorithm and the Gohberg-Semencul-type formula need to solve $T_m x={\bf e}_1$, where $T_m$ is a Toeplitz matrix and ${\bf e}_1=[1,0,...,0]^\intercal$ is a unit vector. Accordingly, the preconditioned conjugate gradient (PCG) method with Strang's circulant preconditioner \cite{chan-1996} is mentioned to solve the above equation. We set a zero initial guess and a stopping relative residual norm tolerance $tol=10^{-15}$.

For the error of the numerical solution in the spacial direction, we use the solution obtained by the two-level Strang splitting method in the fine grid as the reference solution, signifying ${\bf u}(h_{ref},\tau)$. The $L$-infinite ($L^\infty$) norm error and convergence order in space are defined by
$$
E_s (h,\tau)=
\|{\bf u}(h_{ref},\tau)-{\bf u}(h,\tau)\|_\infty, \qquad
Order_h=
\log_2\Bigg(\frac{E_s (2h,\tau)}{E_s (h,\tau)}\Bigg),
$$
where ${\bf u}(h,\tau)$ is the numerical solution with time grid size $\tau$ and space step size $h$ in the coarse grid.

We denote that ${\bf u}(h,\tau)$ is a numerical solution which is gotten by two-level Strang splitting method, Strang splitting method, or splitting-ADI method. Suppose that the $L^\infty$ norm errors and convergence rates in time are computed by\cite{zhao-2014} :
\begin{equation*}
E_t (h,\tau)=
\|{\bf u}(h,2\tau)-{\bf u}(h,\tau)\|_\infty, \qquad
Order_\tau=
\log_2\Bigg(\frac{E_t (h,2\tau)}{E_t (h,\tau)}\Bigg).
\end{equation*}
\begin{example}
We consider a 2D SFAC equation (\ref{FAC}) with $\varepsilon=0.1$, $T=1$, and $\Omega=(0,2)^2$. The initial condition is given by
\begin{equation*}
u^0(x)=\frac{1}{2}e^{-100\big((x^{(1)}-\frac{2}{3})^2+(x^{(2)}-\frac{2}{3})^2\big)}
+\frac{1}{2}e^{-100\big((x^{(1)}-\frac{4}{3})^2+(x^{(2)}-\frac{4}{3})^2\big)}.
\end{equation*}
\end{example}
First, we verify that the proposed method (\ref{finalstrang}) is second-order convergent in space. Setting $\tau=T$, and taking the numerical solution in (\ref{finalstrang}) with $h_1=h_2=\frac{1}{2048}$ as the benchmark solution. The maximum norms of numerical errors and the corresponding convergence orders are presented in Table 1, and it is clear that the value of space convergence order is approximately two.

Then we test the convergence rate in time. Fixing $m=256$ ($h=\frac{1}{128}$), Table 2 shows the $L^\infty$ norm of numerical error, the convergence rate in time and corresponding CPU times (in units of seconds). It is apparent that all methods are second-order convergence in time and the two-level Strang splitting method takes the least amount of CPU times.
\begin{table}[!htbp]
\begin{center}
\caption{Convergence order of the space in the maximum norm for Example 1.}
\def\temptablewidth{1\textwidth}
{\rule{\temptablewidth}{1pt}}
{\footnotesize
\begin{tabular*}{\temptablewidth}{@{\extracolsep{\fill}}ccccccccc}
&  \multicolumn{2}{c}{$(1.1,1.2)$}
&  \multicolumn{2}{c}{$(1.5,1.5)$}
&  \multicolumn{2}{c}{$(1.2,1.8)$}
&  \multicolumn{2}{c}{$(1.6,1.9)$}\\
\cmidrule(l){2-3}\cmidrule(l){4-5}\cmidrule(l){6-7}\cmidrule(l){8-9}
$h$ & Error & Order & Error & Order & Error & Order & Error & Order \\
\specialrule{0.5pt}{0.5pt}{2pt}
$2^{-4}$ & $2.4190e{-2}$  & -    & $3.0865e{-2}$ & -    & $2.1425e{-2}$ & -    & $2.0645e{-2}$ & - \\
$2^{-5}$ & $7.0129e{-3}$ & 1.79 & $8.9022e{-3}$ & 1.79 & $5.8924e{-3}$ & 1.86 & $5.5039e{-3}$ & 1.91 \\
$2^{-6}$ & $1.8595e{-3}$ & 1.92 & $2.3431e{-3}$ & 1.93 & $1.5245e{-3}$ & 1.95 & $1.4050e{-3}$ & 1.97  \\
$2^{-7}$ & $4.7539e{-4}$ & 1.97 & $5.9622e{-4}$ & 1.97 & $3.8574e{-4}$ & 1.98 & $3.5346e{-4}$ & 1.99 \\
$2^{-8}$ & $1.1875e{-4}$ & 2.00 & $1.4862e{-4}$ & 2.00  & $9.5950e{-5}$ & 2.01 & $8.7705e{-5}$ & 2.01 \\
\end{tabular*} {\rule{\temptablewidth}{1pt}}
}
\end{center}
\end{table}

\begin{table}
\begin{center}
\caption{The numerical results of the two-level Strang splitting method, Strang splitting method, and splitting-ADI method for Example 1.}
\def\temptablewidth{1\textwidth}
{\rule{\temptablewidth}{1pt}}
\setlength{\tabcolsep}{1mm}
{\footnotesize
\begin{tabular*}{\temptablewidth}{@{\extracolsep{\fill}}ccccccccccc}
& 
& \multicolumn{3}{c}{two-level Strang splitting}
& \multicolumn{3}{c}{Strang splitting}
& \multicolumn{3}{c}{splitting-ADI } \\
\cmidrule(l){3-5}\cmidrule(l){6-8}\cmidrule(l){9-11}
$(\alpha_1,{\alpha_2})$ & $\tau$ & Error & Order & CPU & Error & Order & CPU & Error & Order & CPU \\
\specialrule{0.5pt}{0.5pt}{2pt}
		  & 1/100  & $1.2824e{-7}$ & - & 0.48  & $7.9046e{-8}$
		  & -    & 18.24 & $2.7095e{-7}$ & - & 3.54\\
		  & 1/200  & $3.2071e{-8}$ & 2.00 & 0.90 & $1.9762e{-8}$
		  & 2.00 & 36.28 & $6.7737e{-8}$ & 2.00 & 7.07 \\
(1.2,1.4) & 1/400  & $8.0182e{-9}$ & 2.00 & 1.83 & $4.9399e{-9}$
		  & 2.00 & 72.32 & $1.6934e{-8}$ & 2.00 & 14.28 \\
		  & 1/800  & $2.0046e{-9}$ & 2.00 & 3.49 & $1.2355e{-9}$
		  & 2.00 & 147.45 & $4.2322e{-9}$ & 2.00 & 28.45 \\
		  & 1/1600  & $5.0115e{-10}$ & 2.00 & 6.73 & $3.0945e{-10}$
		  & 2.00 & 292.58 & $1.0576e{-9}$ & 2.00 & 56.76 \\
\specialrule{0.5pt}{0.5pt}{2pt}
		  & 1/100  & $3.8688e{-7}$ & - & 0.47 & $3.8691e{-7}$
		  & -    & 17.56 & $1.0447e{-6}$ & - & 3.67\\
		  & 1/200  & $9.6720e{-8}$ & 2.00 & 0.91 & $9.6727e{-8}$
		  & 2.00 & 36.75 & $2.6117e{-7}$ & 2.00 & 7.20 \\
(1.5,1.5) & 1/400  & $2.4180e{-8}$ & 2.00 & 1.89 & $2.4182e{-8}$
		  & 2.00 & 73.31 & $6.5292e{-8}$ & 2.00 & 14.45 \\
		  & 1/800  & $6.0450e{-9}$ & 2.00 & 3.45 & $6.0469e{-9}$
		  & 2.00 & 147.63 & $1.6322e{-8}$ & 2.00 & 30.10 \\
		  & 1/1600  & $1.5109e{-9}$ & 2.00 & 6.79 &$1.4963e{-9}$
		  & 2.01 & 294.76 & $4.0800e{-9}$ & 2.00 & 57.71 \\
\specialrule{0.5pt}{0.5pt}{2pt}
		  & 1/100  & $3.2386e{-6}$ & - & 0.46  & $7.9992e{-7}$
		  & -    & 18.30 & $1.9878e{-6}$ & - & 3.66\\
		  & 1/200  & $8.3546e{-7}$ & 1.95 & 0.89  & $1.9998e{-7}$
		  & 2.00  & 35.91 & $4.9696e{-7}$ & 2.00 & 7.12\\
		  & 1/400  & $2.1056e{-7}$ & 1.99 & 1.79 & $4.9996e{-8}$
		  & 2.00 & 73.31 & $1.2424e{-7}$ & 2.00 & 14.43 \\
(1.2,1.8) & 1/800  & $5.2747e{-8}$ & 2.00 & 3.58 & $1.2498e{-8}$
		  & 2.00 & 147.74 & $3.1059e{-8}$ & 2.00 & 28.91 \\
		  & 1/160  & $1.3193e{-8}$ & 2.00 & 7.11 & $3.1172e{-9}$
		  & 2.00 & 294.46 & $7.7632e{-9}$ & 2.00 & 57.75 \\
\specialrule{0.5pt}{0.5pt}{2pt}
		  & 1/100  & $4.0334e{-6}$ & - & 0.47  & $9.6228e{-7}$
		  & -    & 18.56 & $2.8874e{-6}$ & - & 3.67\\
		  & 1/200  & $1.1152e{-6}$ & 1.85 & 0.83  & $2.4058e{-7}$
		  & 2.00   & 35.92 & $7.2184e{-7}$ & 2.00 & 7.15 \\
		  & 1/400  & $2.8653e{-7}$ & 1.96 & 1.68 & $6.0143e{-8}$
		  & 2.00 & 73.38 & $1.8046e{-7}$ & 2.00 & 14.34\\
(1.6,1.9) & 1/800  & $7.2136e{-8}$ & 1.99 & 3.40 & $1.5034e{-8}$
		 & 2.00 & 155.08 & $4.5115e{-8}$ & 2.00 & 28.57\\
		  & 1/1600  & $1.8066e{-8}$ & 2.00 & 6.88 & $3.7612e{-9}$
		  & 2.00 & 296.49 & $1.1279e{-8}$ & 2.00 & 57.36 \\
\end{tabular*} {\rule{\temptablewidth}{1pt}}
}
\end{center}
\end{table}

\begin{example}
In this example, we consider the SFAC equation in 3D case. We set $\varepsilon=0.1$, $\Omega=(0,1)^3$, and $T=1$. The initial condition is
\begin{equation*}
u^0(x)=\frac{1}{2}e^{-500\big((x^{(1)}-\frac{3}{8})^2+(x^{(2)}-\frac{3}{8})^2+(x^{(3)}-\frac{3}{8})^2\big)}
+\frac{1}{2}e^{-500\big((x^{(1)}-\frac{5}{8})^2+(x^{(2)}-\frac{5}{8})^2+(x^{(3)}-\frac{5}{8})^2\big)}.
\end{equation*}
\end{example}
Firstly, the prerequisites will be presented similarly to Example 1, we fix $\tau=T$ and test the convergence order of space in the presence of different space step size $h$. For the benchmark, we take the proposed method with $h_1=h_2=h_3=\frac{1}{512}$. Table 3 shows the numerical results of the maximum norm of error and convergence rate. In the same situation as Example 1, the spatial convergence rate is approximately two.

Then for testing the time convergence order, refer to Example 1, we select $m=64$ ($h=\frac{1}{64}$) and compare the numerical solution of the two-level Strang splitting method, Strang splitting method, splitting-ADI  and proposed method. Again, from Table 4, we observe that these methods process second-order accurately in time. Moreover, the two-level Strang splitting method is the fastest with concerning the computational cost and has the lowest error.

\begin{table}[!htbp]
\begin{center}
\caption{Convergence orders of the space in the maximum norm for Example 2.}
\def\temptablewidth{1\textwidth}
{\rule{\temptablewidth}{1pt}}
{\footnotesize
\begin{tabular*}{\temptablewidth}{@{\extracolsep{\fill}}ccccccccc}
&  \multicolumn{2}{c}{$(1.2,1.3,1.4)$}
&  \multicolumn{2}{c}{$(1.5,1.5,1.5)$}
&  \multicolumn{2}{c}{$(1.6,1.7,1.8)$}
&  \multicolumn{2}{c}{$(1.2,1.5,1.8)$}\\
\cmidrule(l){2-3}\cmidrule(l){4-5}\cmidrule(l){6-7}\cmidrule(l){8-9}
$h$ & Error & Order & Error & Order & Error & Order & Error & Order \\
\specialrule{0.5pt}{0.5pt}{2pt}
$2^{-4}$ & $9.8378e{-2}$  & -    & $5.1707e{-2}$ & -    & $3.2576e{-2}$ & -    & $1.2156e{-2}$ & - \\
$2^{-5}$ & $2.6424e{-2}$ & 1.90 & $9.1509e{-3}$ & 2.51 & $7.2734e{-3}$ & 2.61 & $1.7605e{-3}$ & 2.79 \\
$2^{-6}$ & $6.4887e{-3}$ & 2.03 & $1.9936e{-3}$ & 2.18 & $1.7809e{-3}$ & 2.03 & $ 4.0192e{-4}$ & 2.13  \\
$2^{-7}$ & $1.5580e{-3}$ & 2.06 & $4.6617e{-4}$ & 2.10 & $4.2834e{-4}$ & 2.06 & $9.4496e{-5}$ & 2.09 \\
\end{tabular*} {\rule{\temptablewidth}{1pt}}
}
\end{center}
\end{table}

\begin{table}[!htbp]
\begin{center}
\caption{Convergence order of the space in the maximum norm for Example 2.}
\def\temptablewidth{1\textwidth}
\setlength{\tabcolsep}{1mm}
{\rule{\temptablewidth}{1pt}}
{\footnotesize
\begin{tabular*}{\temptablewidth}{@{\extracolsep{\fill}}ccccccccccc}
& 
& \multicolumn{3}{c}{two-level Strang splitting}
& \multicolumn{3}{c}{Strang splitting}
& \multicolumn{3}{c}{splitting-ADI } \\
\cmidrule(l){3-5}\cmidrule(l){6-8}\cmidrule(l){9-11}
$(\alpha_1,{\alpha_2},{\alpha_3})$ & $\tau$ & Error & Order & CPU & Error & Order & CPU & Error & Order & CPU \\
\specialrule{0.5pt}{0.5pt}{2pt}
			  & 1/20  & $1.9858e{-6}$ & - & 0.30  & $1.3388e{-5}$
			  & -    & 41.25 & $5.9441e{-5}$ & - & 3.64\\
			  & 1/40  & $4.9648e{-7}$ & 2.00 & 0.59 & $3.3481e{-6}$
			  & 2.00 & 86.72 & $1.4858e{-5}$ & 2.00 & 7.27 \\
(1.2,1.3,1.4) & 1/80  & $1.2412e{-7}$ & 2.00 & 1.20 & $8.3710e{-7}$
			  & 2.00 & 176.87 & $3.7143e{-6}$ & 2.00 & 14.84 \\
			  & 1/160  & $3.1031e{-8}$ & 2.00 & 2.36 & $2.0928e{-7}$
			  & 2.00 & 347.84 & $9.2856e{-7}$ & 2.00 & 29.52 \\
			  & 1/320  & $7.7577e{-9}$ & 2.00 & 4.75 & $5.2320e{-8}$
			  & 2.00 & 735.42 & $2.3214e{-7}$ & 2.00 & 58.84 \\
\specialrule{0.5pt}{0.5pt}{2pt}
			  & 1/20  & $1.2394e{-6}$ & - & 0.28  & $9.3871e{-6}$
			  & -    & 44.31 & $6.5795e{-5}$ & - & 3.66\\
			  & 1/40  & $3.1052e{-7}$ & 2.00 & 0.54 & $2.3530e{-6}$
			  & 2.00 & 91.35 & $1.6458e{-5}$ & 2.00 & 7.16 \\
(1.5,1.5,1.5) & 1/80  & $7.7671e{-8}$ & 2.00 & 1.09 & $5.8866e{-7}$
		  	  & 2.00 & 187.10 & $4.1149e{-6}$ & 2.00 & 14.76 \\
			  & 1/160  & $1.9420e{-8}$ & 2.00 & 2.12 & $1.4719e{-7}$
			  & 2.00 & 368.19 & $1.0287e{-6}$ & 2.00 & 29.65 \\
			  & 1/320  & $4.8553e{-9}$ & 2.00 & 4.28 & $3.6799e{-8}$
			  & 2.00 & 690.41 & $2.5718e{-7}$ & 2.00 & 58.73 \\
\specialrule{0.5pt}{0.5pt}{2pt}
		      & 1/20  & $1.1917e{-5}$ & - & 0.28  & $1.5444e{-5}$
			  & -    & 42.04 & $6.2305e{-5}$ & - & 3.76\\
			  & 1/40  & $3.1912e{-6}$ & 1.90 & 0.54 & $3.9056e{-6}$
			  & 1.99 & 86.39 & $1.5606e{-5}$ & 2.00 & 7.28 \\
(1.2,1.5,1.8) & 1/80  & $8.1238e{-7}$ & 1.97 & 1.06 & $9.7926e{-7}$
			  & 2.00 & 178.17 & $3.9030e{-6}$ & 2.00 & 14.59 \\
			  & 1/160  & $2.0403e{-7}$ & 1.99 & 2.22 & $2.4499e{-7}$
			  & 2.00 & 353.24 & $9.7637e{-7}$ & 2.00 & 28.89 \\
			  & 1/320  & $5.1066e{-8}$ & 2.00 & 4.34 & $6.1259e{-8}$
			  & 2.00 & 712.84 & $2.4405e{-7}$ & 2.00 & 58.48 \\
\specialrule{0.5pt}{0.5pt}{2pt}
			  & 1/20  & $4.9153e{-6}$ & - & 0.29  & $7.3008e{-6}$
			  & -    & 42.29 & $3.9694e{-5}$ & - & 3.59\\
			  & 1/40  & $1.3162e{-6}$ & 1.90 & 0.54 & $1.8636e{-6}$
			  & 1.97 & 85.77 & $9.9616e{-6}$ & 2.00 & 7.31 \\
(1.6,1.7,1.8) & 1/80  & $3.3506e{-7}$ & 1.98 & 1.16 & $4.6845e{-7}$
			  & 1.99 & 175.49 & $2.4935e{-6}$ & 2.00 & 14.57 \\
			  & 1/160  & $8.4152e{-8}$ & 1.99 & 2.17 & $1.1727e{-7}$
			  & 2.00 & 352.74 & $6.2370e{-7}$ & 2.00 & 29.20 \\
			  & 1/320  & $2.1062e{-8}$ & 2.00 & 4.36 & $2.9328e{-8}$
			  & 2.00 & 698.04 & $1.5590e{-7}$ & 2.00 & 57.78 \\
\end{tabular*} {\rule{\temptablewidth}{1pt}}
}
\end{center}
\end{table}

\newpage
\bigskip
\noindent\textbf{Example 3.} In this example, we consider the 2D SFAC equation (\ref{FAC}) with the parameter $\varepsilon=0.01$ on the domain $\Omega=(0,1)^2$ with evenly distributed random starting data between $-0.9$ and $0.9$. We generate on the $255 \times 255$ mesh, and the orders of fractional derivative are set by $(\alpha, {\beta})=(1.1,1.3)$, $(1.5,1.5)$, $(1.6,1.9)$.
For 3D case, we have $127 \times 127 \times 127$ mesh generalization, and choose $(\alpha,{\beta},{\gamma})=(1.1,1.2,1.3)$, $(1.5,1.5,1.5)$, $(1.7,1.8,1.9)$. In both cases, the time step size is $\tau=0.1$. The effects of fractional diffusion on phase separation and coarsening are investigated. The snapshots of the numerical solutions are presented in the Fig. 1 and Fig. 2.

From Fig. 1 and Fig. 2, we observe that reducing the fractional power results in thinner surfaces, smaller bulk areas, and a more heterogeneous phase structure. Moreover, inordinate $\alpha$, the interfaces are thinner than those in abscissa ${\beta}$, and the thinnest axis is in vertical $\gamma$, which means that the phase coarsening process is slower for smaller fractional orders.
\begin{figure}[!htbp]
	\centering
	\begin{minipage}[t]{0.33\linewidth}
		\includegraphics[width=1.5in]{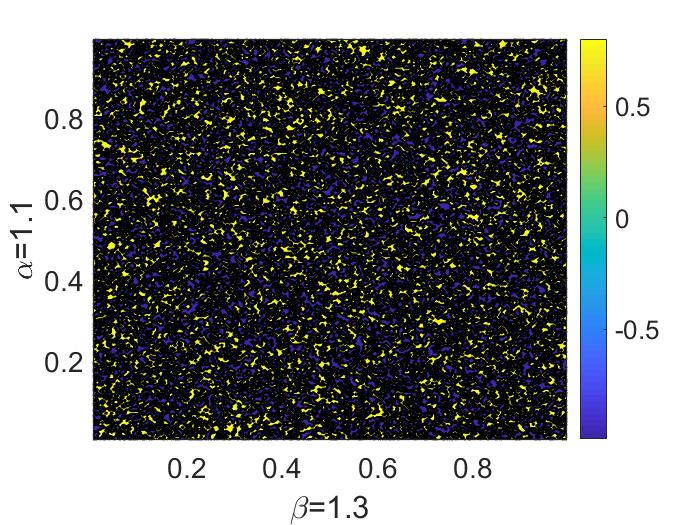}
	\end{minipage}%
	\begin{minipage}[t]{0.33\linewidth}
		\includegraphics[width=1.5in]{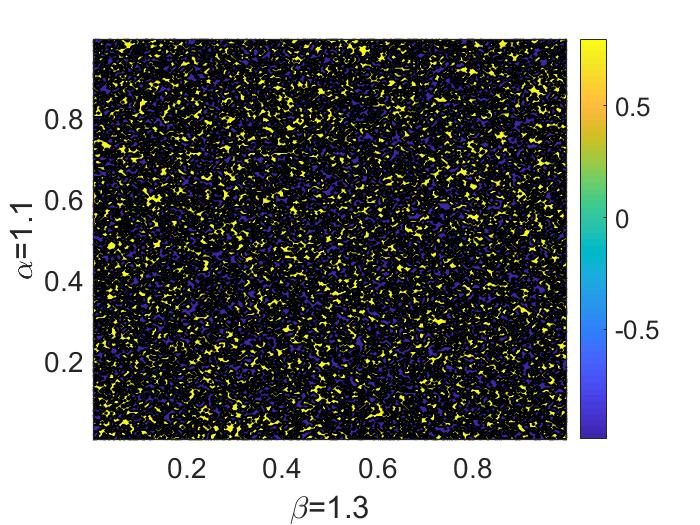}
	\end{minipage}%
	\begin{minipage}[t]{0.33\linewidth}
		\includegraphics[width=1.5in]{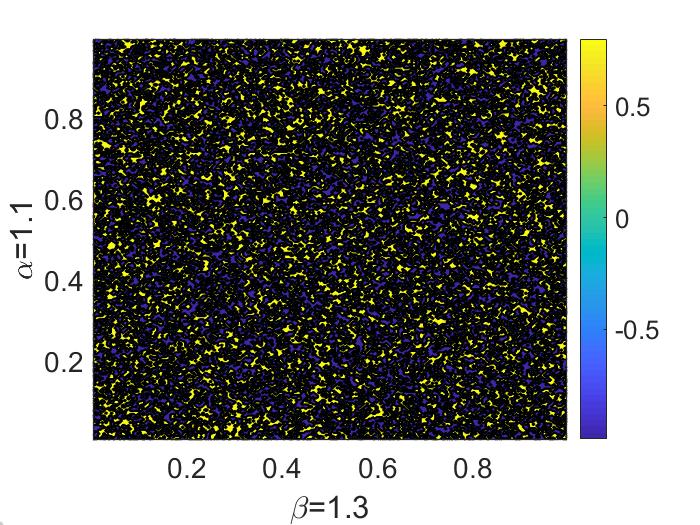}
	\end{minipage}
	\hfill
	\begin{minipage}[t]{0.33\linewidth}
		\includegraphics[width=1.5in]{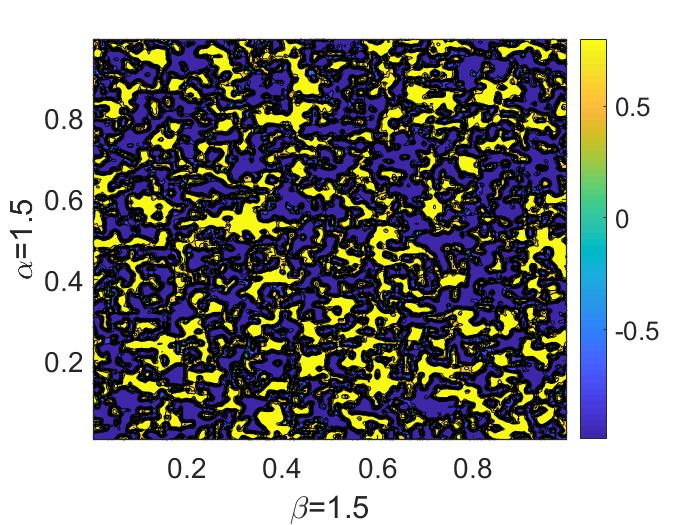}
	\end{minipage}%
	\begin{minipage}[t]{0.33\linewidth}
		\includegraphics[width=1.5in]{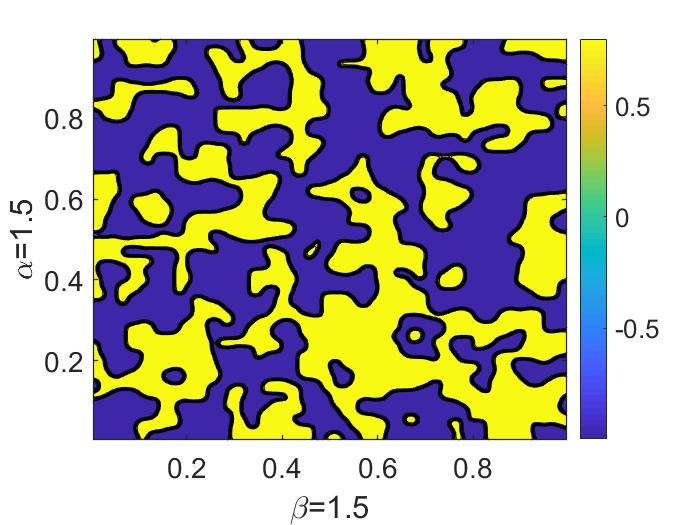}
	\end{minipage}%
	\begin{minipage}[t]{0.33\linewidth}
		\includegraphics[width=1.5in]{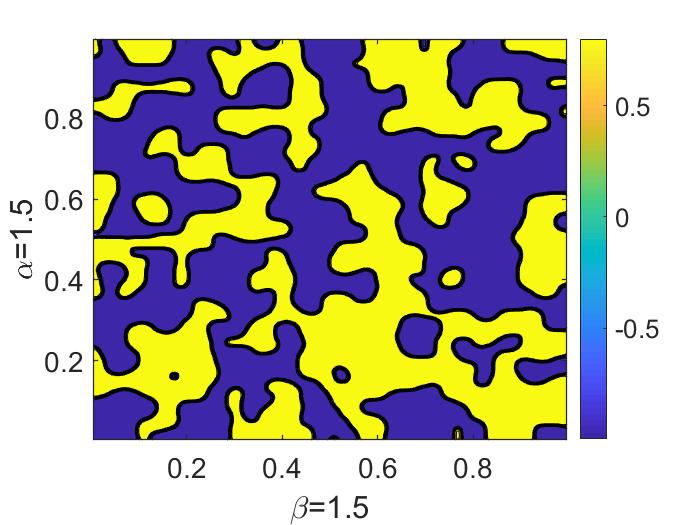}
	\end{minipage}
	\hfill
	\begin{minipage}[t]{0.33\linewidth}
		\includegraphics[width=1.5in]{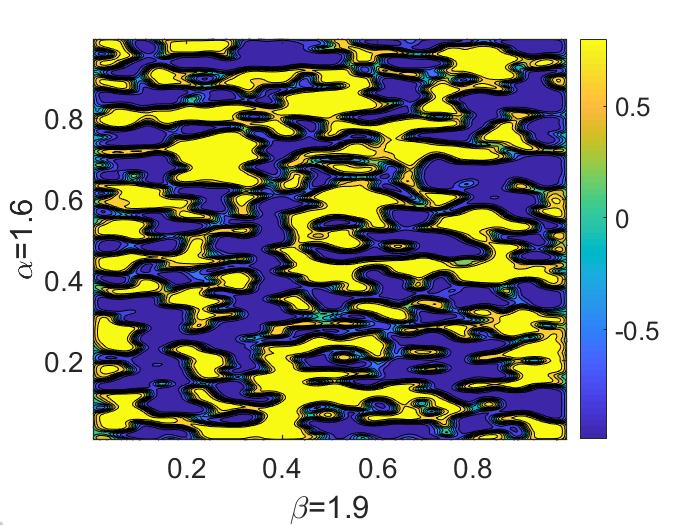}
	\end{minipage}%
	\begin{minipage}[t]{0.33\linewidth}
		\includegraphics[width=1.5in]{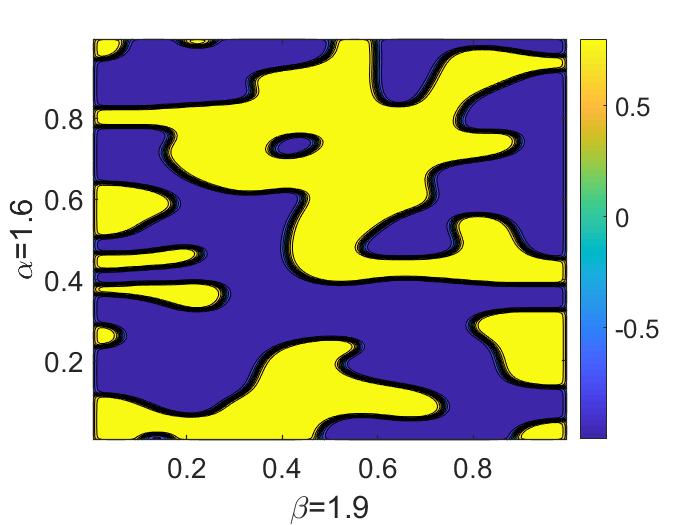}
	\end{minipage}%
	\begin{minipage}[t]{0.33\linewidth}
		\includegraphics[width=1.5in]{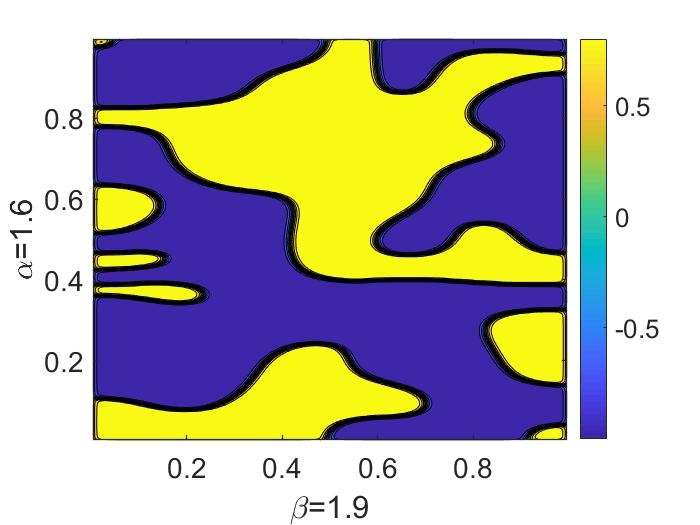}
	\end{minipage}
	\centering
	\caption{Evolution of the numerical solution for Example 3 by the two-level Strang splitting scheme. From left to right: $t=5,60,$ and $100$.}
\end{figure}

\begin{figure}[!htbp]
	\centering
	\begin{minipage}[t]{0.33\linewidth}
		\includegraphics[width=1.5in]{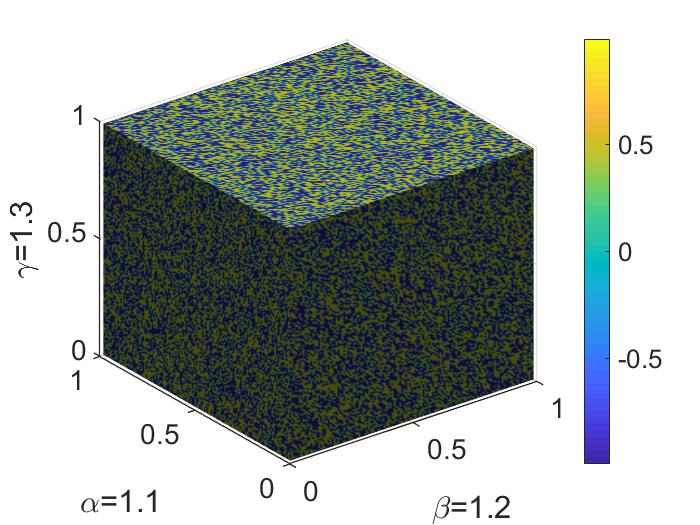}
	\end{minipage}%
	\begin{minipage}[t]{0.33\linewidth}
		\includegraphics[width=1.5in]{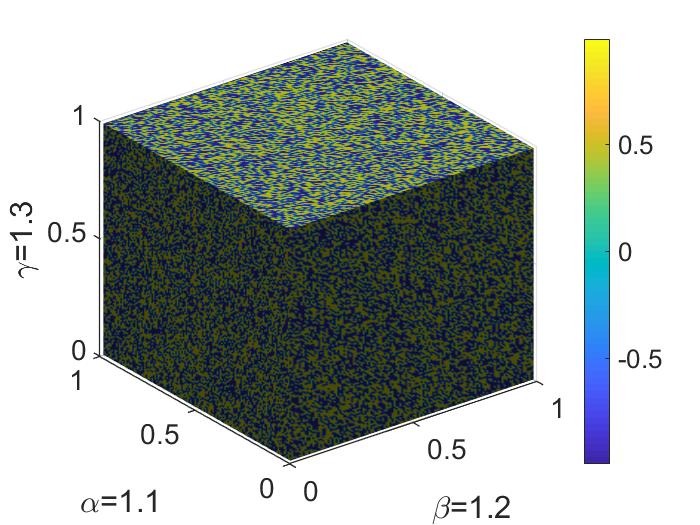}
	\end{minipage}%
	\begin{minipage}[t]{0.33\linewidth}
		\includegraphics[width=1.5in]{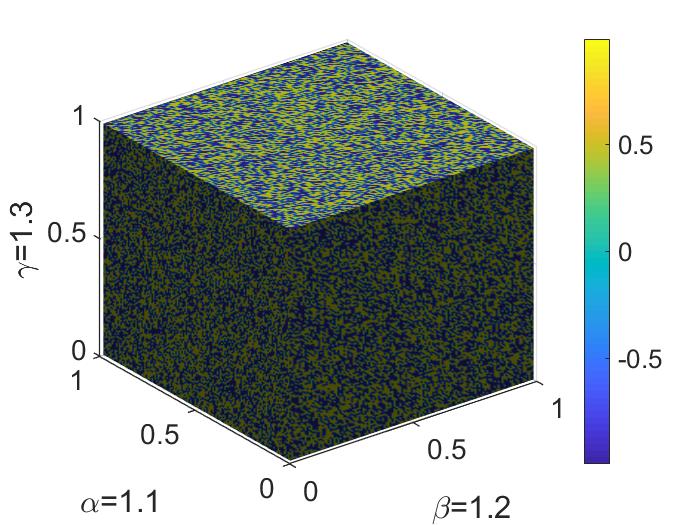}
	\end{minipage}
	\hfill
	\begin{minipage}[t]{0.33\linewidth}
		\includegraphics[width=1.5in]{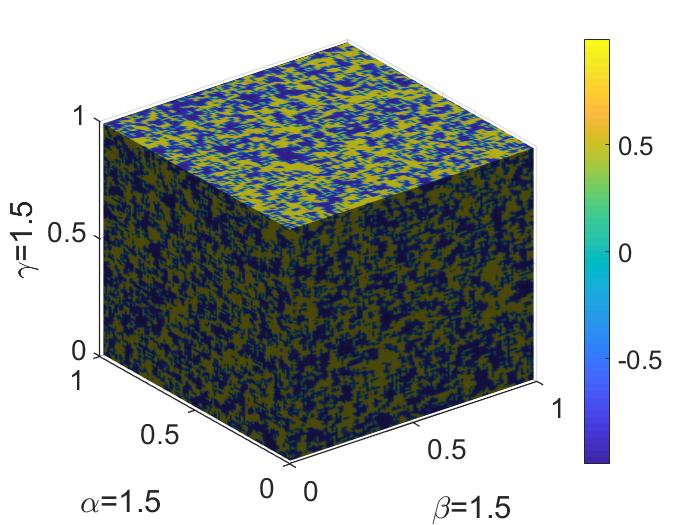}
	\end{minipage}%
	\begin{minipage}[t]{0.33\linewidth}
		\includegraphics[width=1.5in]{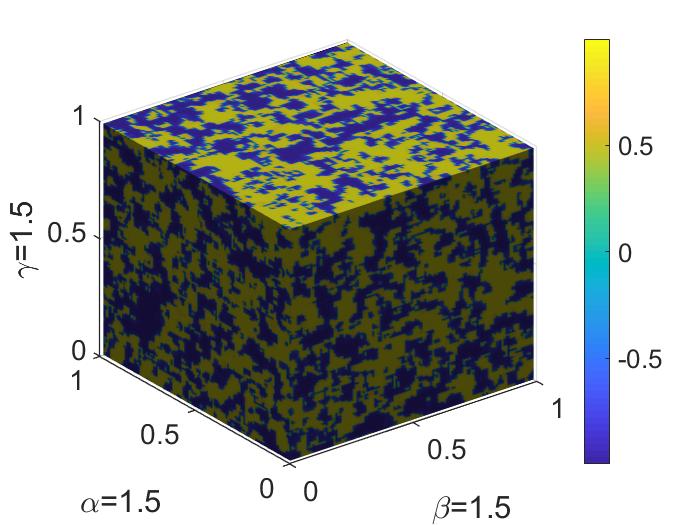}
	\end{minipage}%
	\begin{minipage}[t]{0.33\linewidth}
		\includegraphics[width=1.5in]{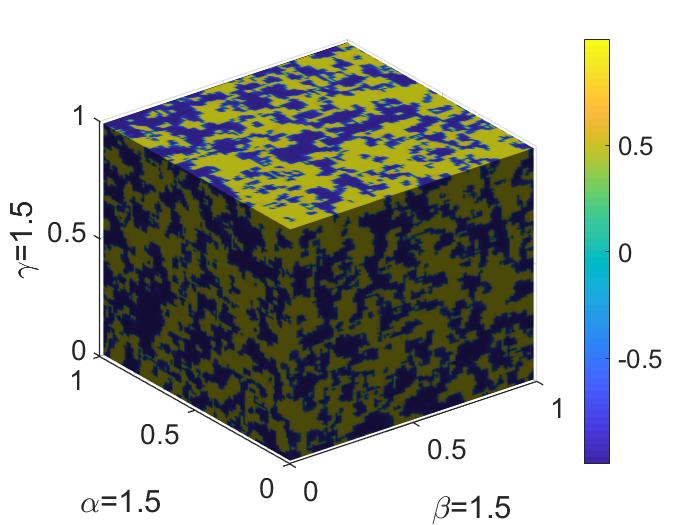}
	\end{minipage}
	\hfill
	\begin{minipage}[t]{0.33\linewidth}
		\includegraphics[width=1.5in]{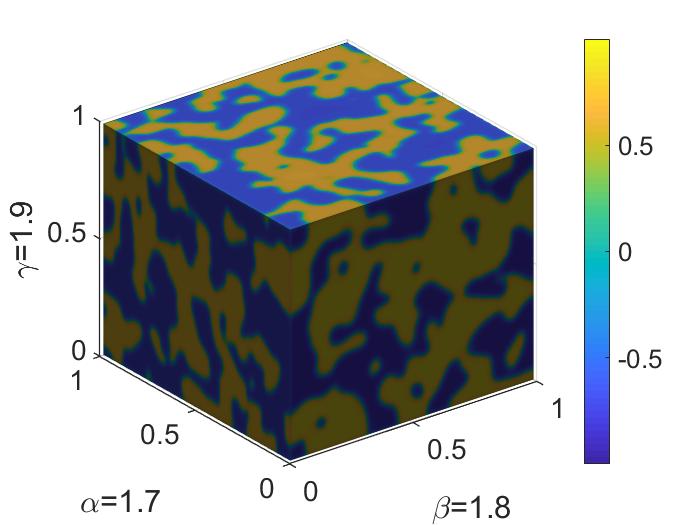}
	\end{minipage}%
	\begin{minipage}[t]{0.33\linewidth}
		\includegraphics[width=1.5in]{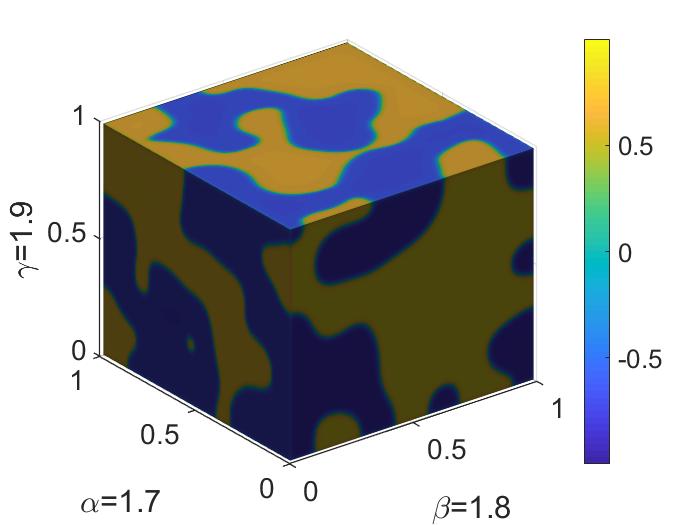}
	\end{minipage}%
	\begin{minipage}[t]{0.33\linewidth}
		\includegraphics[width=1.5in]{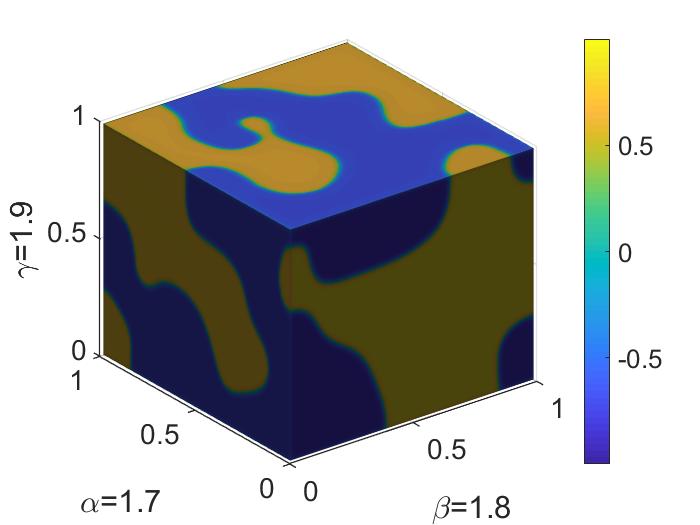}
	\end{minipage}
\caption{Evolution of the numerical solution for Example 3 by the two-level Strang splitting scheme. From left to right: $t=15,100,$ and $200$.}
\end{figure}

\newpage
\bigskip
\noindent\textbf{Example 4.} In this example, we will show the discrete maximum principle of the numerical scheme (\ref{finalstrang}). Consider the 2D/3D SFAC equation (\ref{FAC}) with initial condition
$$
u_{\bf{x}}(0)=0.95 \times rand(x^{(d)}) +0.05,
$$
where $rand(\cdot)$ is a random number in (0,1) and $x^{(d)}$ means the d-dimensional elements. All conditions are the same as the above Example 3, except for the $\varepsilon=0.1$ and the value of the time step is $\tau=0.01$. The maximum norms at each moment are plotted in Fig. 3. Our observations revealed that the scheme (\ref{finalstrang}) satisfies the discrete maximum principle under the condition $\|u^0(x)\|_{\infty}\leq 1$. Also, the maximum norm will rise as the fractional orders increase but not surpass 1.
\begin{figure}[!htbp]
	\centering
	\begin{minipage}[t]{0.33\linewidth}
		\includegraphics[width=1.5in]{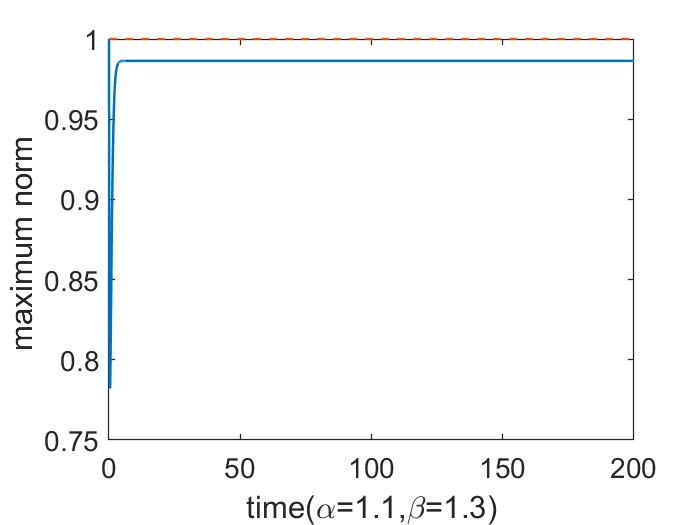}
	\end{minipage}%
	\begin{minipage}[t]{0.33\linewidth}
		\includegraphics[width=1.5in]{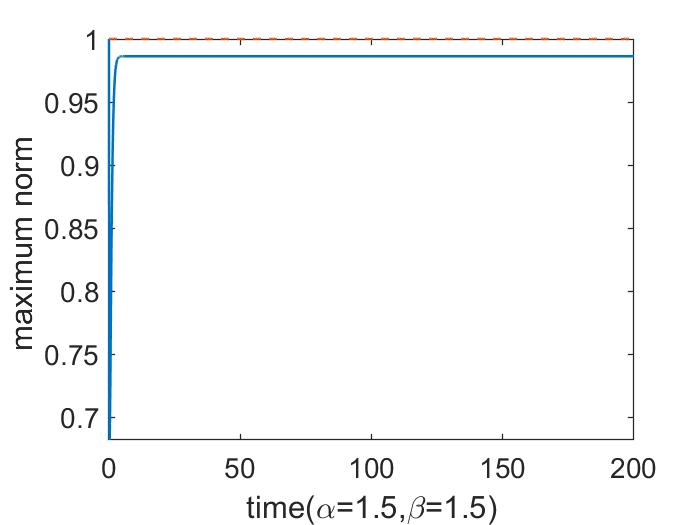}
	\end{minipage}%
	\begin{minipage}[t]{0.33\linewidth}
		\includegraphics[width=1.5in]{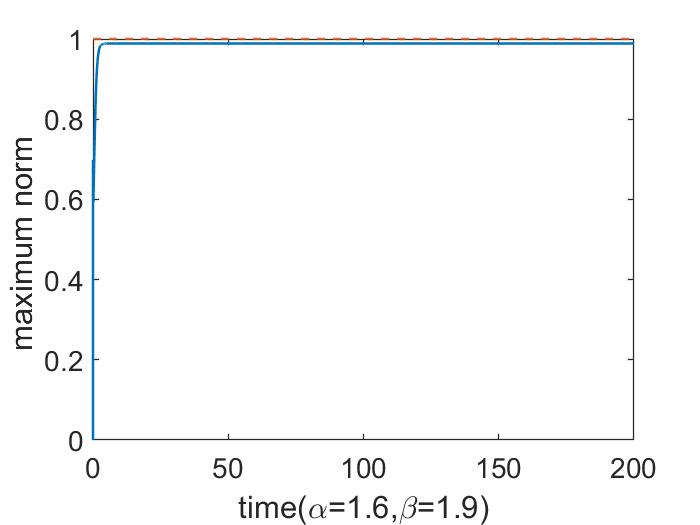}
	\end{minipage}
	\hfill
	\begin{minipage}[t]{0.33\linewidth}
		\includegraphics[width=1.5in]{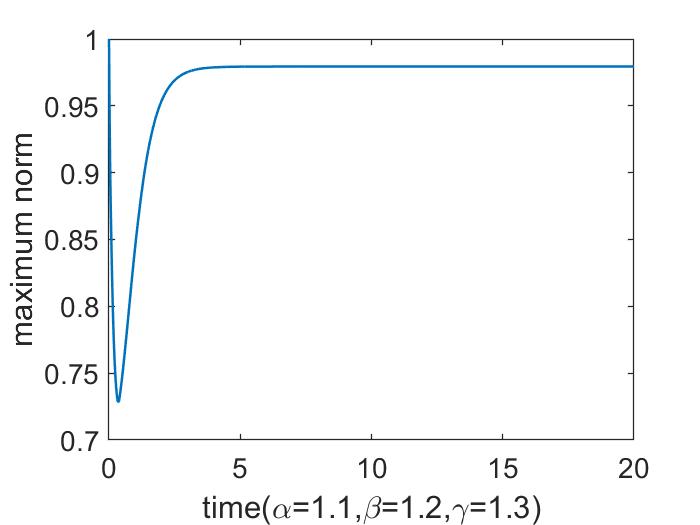}
	\end{minipage}%
	\begin{minipage}[t]{0.33\linewidth}
		\includegraphics[width=1.5in]{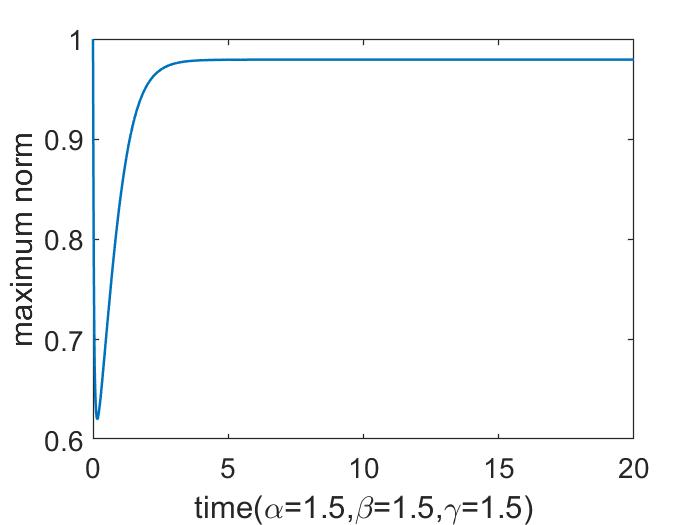}
	\end{minipage}%
	\begin{minipage}[t]{0.33\linewidth}
		\includegraphics[width=1.5in]{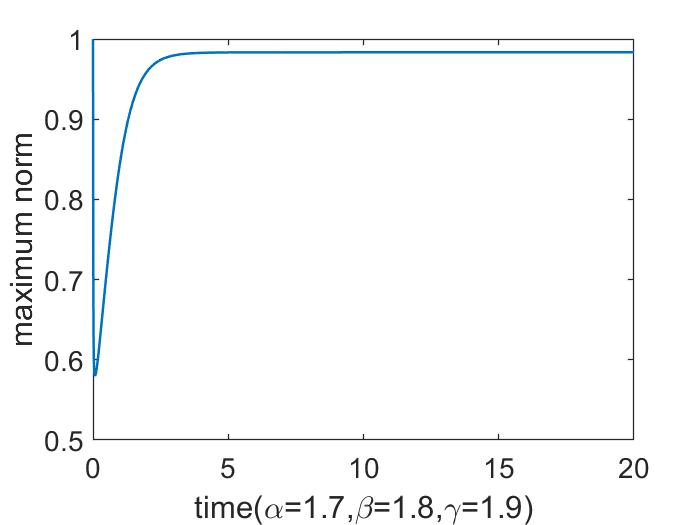}
	\end{minipage}
\caption{Evolution of the maximum norm for Example 4 by the two-level Strang splitting scheme. The 2D situation is displayed by the top three pictures, while the bottom three pictures are the situations of the 3D case.}
\end{figure}

\bigskip
\noindent\textbf{Example 5.}
In this example, we consider the 2D SFAC equation (\ref{FAC}) on the domain $\Omega=(0,1)^2$ with the initial condition
$$
u_{\bf{x}}(0)=0.8 \times rand(x^{(1)},x^{(2)})-0.4.
$$
The goal of this study is to see how the magnitude of the time step affects the energy norm and maximum norm. We fix $\tau=0.01$,  $h_1=h_2=\frac{1}{128}$, $\varepsilon=0.01$, and the fractional order $(\alpha,{\beta},{\gamma})$ are the same as the previous examples. For the 3D case, we set
$$
u_{\bf{x}}(0)=1 \times rand(x^{(1)},x^{(2)},x^{(3)})-0.5.
$$
Here we choose $\varepsilon=0.01$, the uniform mesh $127 \times 127 \times 127$ with the time step $\tau=0.01$.

The fractional Allen-Cahn equation can be thought of as the $L^2$-gradient flow of the fractional Ginzburg-Landau free energy functional
\begin{equation}\label{energy}
E(u)=\int_{\Omega}\big(g(u)-\frac{1}{2}\varepsilon^2 u\mathcal{L}_{x^d}^\alpha u\big)du,
\end{equation}
where $g(u)=\frac{1}{4}(u^2-1)^2$. The discrete energy is represented by
$$
E\big({\bf u}(t_n)\big)
=h^2\Big(\frac{1}{4}\big({\bf u}(t_n)^{2}-1\big)^\intercal
\big({\bf  u}(t_n)^{2}-1\big)
-\frac{1}{2} {\bf u}(t_n)^\intercal A{\bf u}(t_n)\Big),
$$
and we ignore the formula of the 3D case due to analogy. Fig.4 shows the energy decays monotonically.

\begin{figure}[!htbp]
	\begin{minipage}[t]{0.33\linewidth}
		\centering
		\includegraphics[width=1.5in]{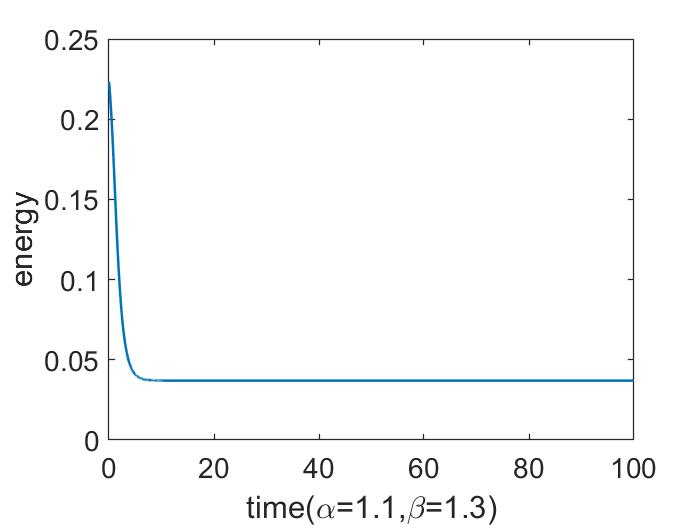}
	\end{minipage}%
	\begin{minipage}[t]{0.33\linewidth}
		\centering
		\includegraphics[width=1.5in]{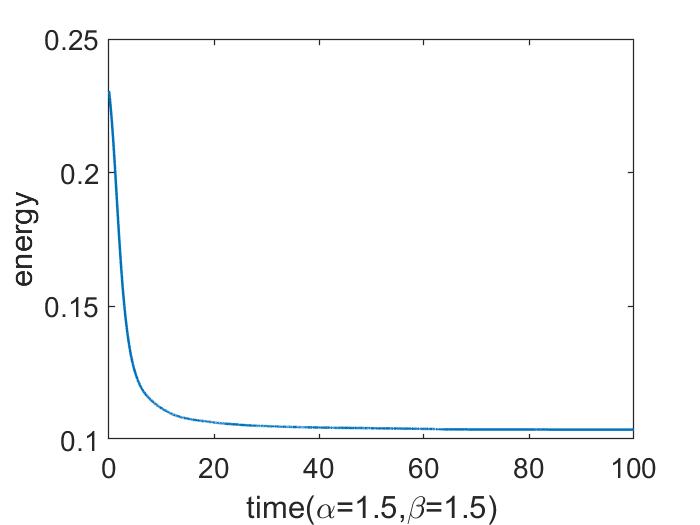}
	\end{minipage}%
	\begin{minipage}[t]{0.33\linewidth}
		\centering
		\includegraphics[width=1.5in]{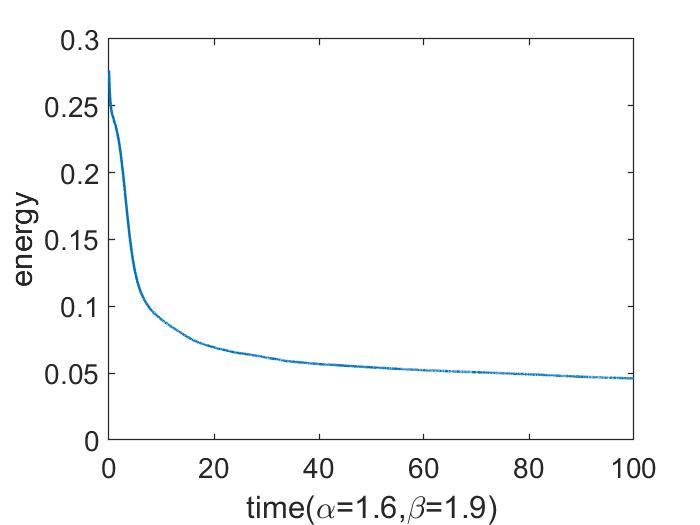}
	\end{minipage}
	\hfill
	\begin{minipage}[t]{0.33\linewidth}
		\centering
		\includegraphics[width=1.5in]{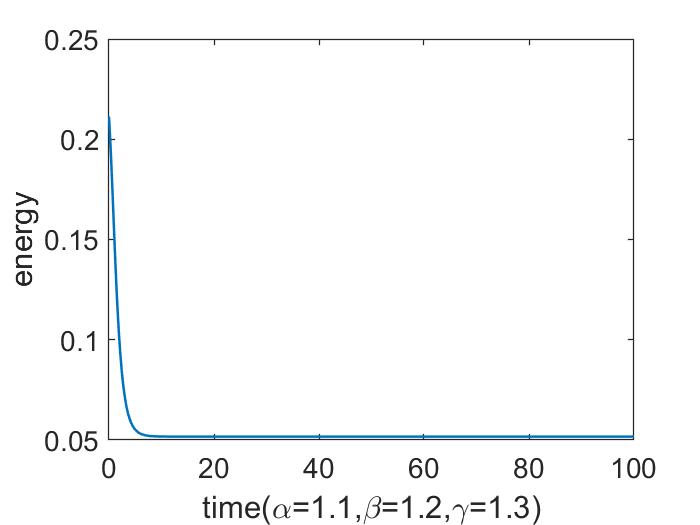}
	\end{minipage}%
	\begin{minipage}[t]{0.33\linewidth}
		\centering
		\includegraphics[width=1.5in]{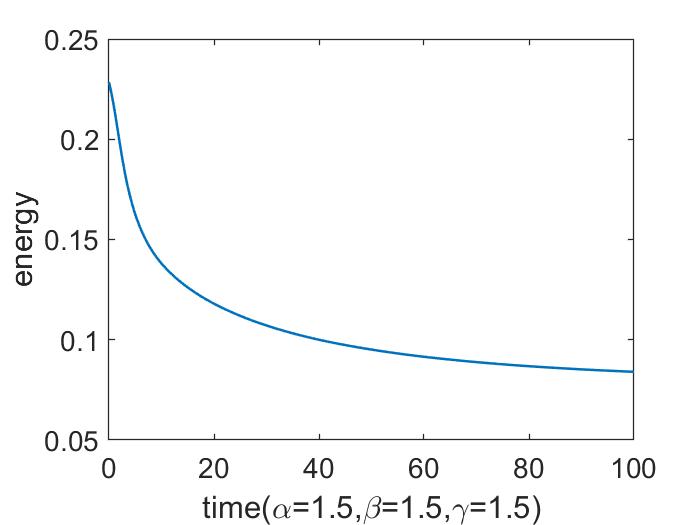}
	\end{minipage}%
	\begin{minipage}[t]{0.33\linewidth}
		\centering
		\includegraphics[width=1.5in]{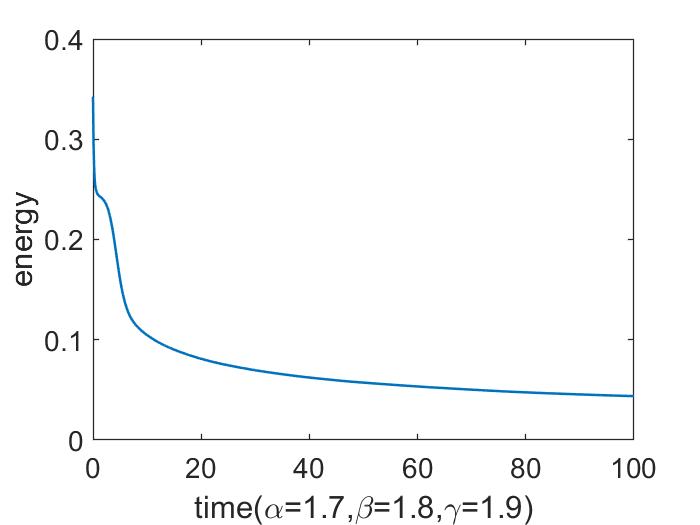}
	\end{minipage}
\caption{Evolution of the energy for Example 4 by the two-level Strang splitting scheme. The 2D situation is displayed by the top three pictures, while the bottom three pictures are the situations of the 3D case.}
\end{figure}

\section{Concluding remarks}\label{S-conclude}
In this paper, we have employed a two-level Strang splitting method for the SFAC equation (\ref{FAC}). After a second-order finite difference for space discretization, the resulting semi-discretized system is solved by the two-level Strang splitting algorithm where the linear subproblem in original splitting method is derived into the circulant and skew-circulant matrices based Strang splitting solver. The significance of this work is reducing the computational cost for solving the SFAC equation by FFTs. Theoretically, we prove that the proposed method satisfies the discrete maximum principle and the convergence order possesses the second order in space and time. Numerical tests have verified the theoretical proofs and shown some properties of physics.

In our future consideration, it is interesting to extend the two-level method to other fractional operators in phase-field equations and check if there are any additional improvements through this method.

\section*{Acknowledgments}
\noindent Funding: This work is supported in part by research grants of the Basic and Applied Basic Research Foundation of Guangdong Province (file no. 2019A1515110893), the National Natural Science Foundation (file no. 12101123), the Science and Technology Development Fund, Macau SAR (file no. 0122/2020/A3), and University of Macau (file no. MYRG2020-00224-FST).


\end{document}